\documentclass{siamltex}
\usepackage{animate}
\usepackage{amsmath}
\usepackage{amssymb}
\usepackage{url}
\usepackage{graphics}
\usepackage{subfigure}
\usepackage{graphicx}
\usepackage{textcomp}
\usepackage{mathrsfs}
\usepackage{epstopdf}
\usepackage{algorithmic, algorithm}
\usepackage{braket,amsfonts,amsopn}
\usepackage{color}
\usepackage{mathscinet}
\usepackage{enumitem}
\usepackage{tikz}
\usepackage{cite}
\usepackage[T1]{fontenc}
\usepackage{lmodern}
\usepackage{dsfont,bbm}

\usepackage{hyperref}
\usepackage{cleveref}
\usepackage{catoptions}
 \AtBeginDocument{%
    \crefname{equation}{equation}{equations}%
    \crefname{chapter}{chapter}{chapters}%
    \crefname{section}{section}{sections}%
    \crefname{appendix}{appendix}{appendices}%
    \crefname{enumi}{item}{items}%
    \crefname{footnote}{footnote}{footnotes}%
    \crefname{figure}{figure}{figures}%
    \crefname{table}{table}{tables}%
    \crefname{theorem}{theorem}{theorems}%
    \crefname{lemma}{lemma}{lemmas}%
    \crefname{corollary}{corollary}{corollaries}%
    \crefname{proposition}{proposition}{propositions}%
    \crefname{definition}{definition}{definitions}%
    \crefname{result}{result}{results}%
    \crefname{example}{example}{examples}%
    \crefname{remark}{remark}{remarks}%
    \crefname{note}{note}{notes}%
}

\numberwithin{equation}{section}
\numberwithin{theorem}{section} % important bit

\newtheorem{rem}{Remark}

\title{Gradient recovery for elliptic interface problem: I. body-fitted mesh}
\author{Hailong Guo%
\thanks{Department of Mathematics, University of California Santa Barbara, CA, 93106 (hlguo@math.ucsb.edu).}%
\and
Xu Yang %
\thanks{Department of Mathematics, University of California Santa Barbara, CA, 93106 (xuyang@math.ucsb.edu).}%
}

\begin{document}

\maketitle

% REQUIRED
\begin{abstract}
In this paper, we propose a novel gradient recovery method for elliptic interface problem
using body-fitted mesh in two dimension. Due to the lack of regularity of solution at interface,
standard gradient recovery methods fail to give superconvergent results, and thus will lead to
overrefinement when served as {\it a posteriori} error estimator. This drawback is overcome by
designing an immersed gradient recovery operator in our method. We prove the superconvergence of this method for both
mildly unstructured mesh and adaptive mesh, and present several numerical examples to verify the superconvergence
and its robustness as {\it a posteriori} error estimator.
\end{abstract}

% REQUIRED

% REQUIRED

\section{Introduction}

Elliptic interface problem frequently appears in the fields of fluid dynamics and material science, where
background consists of rather different materials. The numerical challenge comes from discontinuities of coefficient
at interface, where solution is not smooth in general. Computational methods for elliptic interface problem have been studied intensively in literature, which can be roughly categorized into two types: unfitted mesh methods and body-fitted mesh methods.

Numerical methods based on unfitted mesh solve interface problems on Cartesian grids, among which, famous examples include immersed boundary method (IBM) by Peskin \cite{Peskin1977, Peskin2002} and immersed interface method (IIM) by Leveque and Li \cite{LevequeLi1994}, just to name a few. We refer interested readers to \cite{LiIto2006} for a review of the literature. IBM uses Dirac $\delta$-function to model discontinuity and discretizes it to distribute a singular source to nearest grid point. IIM constructs a special finite difference scheme near interface to get an accurate approximation of the solution. It was further developed in the framework of finite element method \cite{Li1998, LiLinLin2004, LiLinWu2003}, which modifies basis functions on interface elements. Moreover, in \cite{HouLiu2005,HouWuZhang2004}, a special weak form was derived based on Petro-Galerkin method to discretize elliptic interface problem. A shortcoming of unfitted mesh methods is that, the resulting discretized linear system is in general non-symmetric and indefinite even thought the original continuous problem is self-adjoint.

Body-fitted mesh methods require mesh grids to align with interface in order to capture discontinuity.  The resulting discretized
linear system is symmetric and positive definite if the original continuous problem is self-adjoint.
%; see for example \cite{Borgers1990,WeiChenHuangZheng2014, XieItoLi2008, XieLiQiao2011,Zheng2016}.
Error estimates for finite element method with body-fitted mesh have been established by \cite{Babuska1970,BrambleKing1996,ChenZou1998,Xu1982}.
In particular, \cite{ChenZou1998} showed that smooth interface can be approximated by linear interpolation of distinguished points on interface. Although the solution to interface problem has low global regularity,  the finite element approximation was shown to have nearly the same optimal error estimates in both $L^2$ and energy norms as for regular (non-interface) problems.

Meanwhile, superconvergence analysis has attracted considerable attention in the community of finite element method, and theories have been well developed for regular problems \cite{Babuska2001, Chen2001, Wahlbin1995, ZhuLin1989}. Then it is natural to ask if one can obtain similar superconvergence results for elliptic interface problem. However, limited work has been done in this direction due to the lack of regularity of solution at interface. Recently, \cite{Chou2012,Chou2015}  proposed two special interpolation formula to recover flux for linear and quadratic immersed finite element method in one dimension. Supercloseness was established between finite element solution and linear interpolation of the true solution in \cite{WeiChenHuangZheng2014}.

In this paper, we aim to develop gradient recovery methods for elliptic interface problem based on body-fitted finite element discretization.
Standard gradient recovery operators, including superconvergent patch recovery (SPR)  \cite{ZZ1992a, ZZ1992b} and polynomial preserving recovery (PPR) \cite{ZhangNaga2005, NagaZhang2004, NagaZhang2005}, produce superconvergent recovered gradient only when the solution is smooth enough. Therefore, they can not be applied directly to elliptic interface problem since the solution has low regularity at the interface due to the discontinuity of coefficients. Futuremore, building up a recovery-type {\it a posteriori} error estimator based on these methods will lead to overrefine regions as studied in \cite{CaiZhang2009}.

An observation that we rely on is that, even though the solution has low global regularity, it is piecewise smooth on each subdomain separated by the smooth interface. This motivates us to develop a novel gradient recovery method by applying PPR gradient operator on each subdomain since PPR is a local gradient recovery method.  One one hand, for a node away from interface, we use stand PPR gradient recovery operator; On the other hand, for a node close to interface, we design the gradient recovery operator by fitting a quadratic polynomial in least-squares sense only using the sampling points in each subdomain. This will generate two approximations of gradient in each subdomain for a node on interface, which is consistent with the fact that the solution in general is not continuously differentiable at interface. The method is more like to use a divide-and-conquer strategy, which has also been used in the immersed finite element method \cite{Li1998, LiLinLin2004, LiLinWu2003}.

We prove that the proposed gradient recovery method has superconvergence for the following two types of meshes: Benefited from \cite{WeiChenHuangZheng2014} on the approximation estimate and supercloseness, we are able to establish the superconvergence theory on mildly unstructured meshes; Using the practical assumption and supercloseness results in \cite{WuZhang2007} for adaptive mesh, we show that the proposed recovered gradient method is superconvergent to exact gradient on adaptive mesh. Therefore, the method provides an asymptotically exact
{\it a posteriori} error estimator for elliptic interface problem. Compared to the {\it a posteriori} error estimator in \cite{BernardiVerfurth2000,ChenDai2002,MorinNochetto2003} and recovery-type error estimator in \cite{CaiZhang2009, CaiZhang2010}, the estimator based on the proposed gradient recovery is easier in implementation and asymptotically more exact, which will be verified by several two-dimensional numerical examples.

The rest of the paper is organized as follows. In Section~2, we introduce elliptic interface problem and its finite element approximation based on body-fitted mesh. In Section~3, we first give a brief introduction to polynomial preserving recovery method, based on which, we develop a novel gradient recovery method for elliptic interface problem. In Section~4, superconvergence is proved for the proposed gradient recovery operator on both mildly unstructured mesh and adaptive refined mesh. In addition, we show that the method provides an asymptotically exact {\it a posteriori} error estimator for elliptic interface problem.  In Section~5, serval numerical examples are presented to confirm our theoretical results.  Conclusive remarks are made in Section~6.

\section{Finite element method for elliptic interface problem}
In this section, we first introduce elliptic interface problem, and then describe the finite element approximation using body-fitted mesh.

\subsection{Elliptic interface problem}
Let $\Omega$ be a bounded polygonal domain with  Lipschitz boundary
$\partial \Omega$ in $\mathbb{R}^2$.   A $C^2$-curve $\Gamma$ divides $\Omega$ into two disjoint subdomains $\Omega^-$ and $\Omega^+$, which is typically characterized by zero level set of some level set function $\phi$ \cite{Osher2003, Sethian1996}.  Then $\Omega^- = \{z\in \Omega|\phi(z) <0\}$
and $\Omega^+ = \{z\in \Omega|\phi(z) >0\}$.
We shall consider the following elliptic interface problem
\begin{align}
  -\nabla \cdot (\beta(z) \nabla u(z)) &= f(z),  \quad z \text{ in } \Omega\setminus \Gamma, \label{equ:model}\\
   u & = 0, \quad\quad\,\, z \text{ on } \partial\Omega, \label{equ:bnd}
\end{align}
where the diffusion coefficient $\beta(z) \ge \beta_0$ is a piecewise smooth function, i.e.
\begin{equation}
\beta(z) =
\left\{
\begin{array}{ccc}
    \beta^-(z) &  \text{if } z\in \Omega^-, \\
   \beta^+(z)  &   \text{if } z\in \Omega^+,
\end{array}
\right.
\end{equation}
which has a finite jump of function values across the interface $\Gamma$. At the interface $\Gamma$,
one has the following jump conditions
\begin{align}
   [u]_{\Gamma} &= u^+-u^-=0, \label{equ:valuejump}\\
   [\beta u_n]_{\Gamma} &= \beta^+u_n^+ - \beta^-u^-_n = g, \label{equ:fluxjump}
\end{align}
where $u_n$ denotes the normal flux $\nabla u\cdot n$ with $n$ as the unit outer normal vector of the interface $\Gamma$.

{\bf Notations.}  Let  $C$ denote a generic positive constant which may be different at different occurrences.
For the sake of simplicity, we use $x\lesssim y$ to mean that $x\leq Cy$ for some constant $C$
independent of mesh size.  Standard
notations for Sobolev spaces and their associate norms given in \cite{BrennerScott2008, Ciarlet2002, Evans2008}  are adopted in this paper. Moreover, for a subdomain $A$
of $\Omega$, let $\mathbb{P}_m(A)$ be the space of polynomials of
degree less than or equal to $m$ in $A$ and $n_m$ be the
dimension of $\mathbb{P}_m(A)$ which equals to $\frac{1}{2}(m+1)(m+2)$.
$W^{k,p}(A)$ denotes the Sobolev space with norm
$\|\cdot\|_{k, p, A} $ and seminorm $|\cdot|_{k, p,A}$.
 When $p = 2$, $W^{k,2}(A)$ is simply  denoted by $H^{k}(A)$
 and the subscript $p$ is omitted in its associate norm and seminorm.
  As in \cite{WeiChenHuangZheng2014}, denote
 $W^{k,p}(\Omega^-\cup\Omega^+)$  as the function space consisting of piecewise Sobolev  function $w$  such
 that $w|_{\Omega^-}\in W^{k,p}(\Omega^-)$ and $w|_{\Omega^+}\in W^{k,p}(\Omega^+)$.  For the function
 space $W^{k,p}(\Omega^-\cup\Omega^+)$, define norm  as
 \begin{equation*}
\|w\|_{k,p, \Omega^-\cup\Omega^+} = \left( \|w\|_{k,p, \Omega^-}^p + \|w\|_{k,p, \Omega^+}\right)^{1/p},
\end{equation*}
and seminorm as
 \begin{equation*}
|w|_{k,p, \Omega^-\cup\Omega^+} = \left( |w|_{k,p, \Omega^-}^p + |w|_{k,p, \Omega^+}\right)^{1/p}.
\end{equation*}

The variational formulation of elliptic interface problem  \cref{equ:model,equ:bnd,equ:valuejump,equ:fluxjump} is given by
finding $u \in H^1_0(\Omega)$ such that
 \begin{equation}
(\beta\nabla u, \nabla v)  = (f, v) - \langle g, v\rangle, \quad \forall v \in H^1_0(\Omega),
\label{equ:var}
\end{equation}
where $(\cdot, \cdot)$ and  $\langle g, v\rangle$ are standard $L_2$-inner product in the spaces $L^2(\Omega)$ and
$L^2(\Gamma)$ respectively.
By the positiveness of $\beta$, Lax-Milgram  Theorem implies \cref{equ:var} has  a unique solution.
\cite{ChenZou1998, RS1969} proved that $u\in H^{r}(\Omega^-\cup\Omega^+)$ for $0 \le r\le 2$ and
\begin{equation}
\|u\|_{r, \Omega^-\cup\Omega^+} \lesssim \|f\|_{0, \Omega} + \|g\|_{r-3/2, \Gamma},
\end{equation}
if $f\in L^2(\Omega)$ and $g \in H^{r- 3/2}(\Gamma)$.

\begin{rem}  For the sake of easing  theoretical analysis,   we simply assume homogeneous jump of function value.  In fact, one can
extend the method to inhomogeneous jump of  function value  $[u]_{\Gamma} =q$ by defining a piecewise smooth function $\hat{q}$ that satisfies
$\hat{q}|_{\Gamma}=q$ and $\hat{q}|_{\partial\Omega}=0$, and then the problem \cref{equ:model,equ:bnd,equ:fluxjump}
 is equivalent to find   $u =w+\hat{q}$  with $w\in H^1_0(\Omega)$ such that
  \begin{equation*}
(\beta\nabla w, \nabla v)  = (f, v) - \langle g, v\rangle -  (\beta\nabla \hat{q}, \nabla v), \quad \forall v \in H^1_0(\Omega).
\end{equation*}

\end{rem}

\subsection{Finite element approximation}\label{ssec:fem}
 Denote $\mathcal{T}_h$ to be a body-fitted triangulation of $\Omega$, then every triangle $T \in \mathcal{T}_h$ belongs to one of the following
three different cases:
\begin{enumerate}[label=(\alph*)]
\item $T\subset \overline{\Omega^-}$;
\item $T\subset \overline{\Omega^+}$;
\item $T\cap \Omega^- \neq \emptyset$ and $T\cap \Omega^- \neq \emptyset$, then two of vertices of $T$ lie on $\Gamma$.
\end{enumerate}
For any $T\in \mathcal{T}_h$, denote its diameter and supermum of the diameters of the circles inscribed in $T$ by $h_T$ and $\rho_T$ respectively.
Let $h = \max_{T\in \mathcal{T}_h}h_T$.  Assume that the triangulation of $\Omega$ is shape-regular in the sense that there is a constant $\xi$ such
that $ \frac{h_T}{\rho_T} \le \xi$ for all $T\in \mathcal{T}_h$. Denote $\Gamma_h$ as an approximation to $\Gamma$ which consists of the edges with both endpoints lying on $\Gamma$. The domain $\Omega$ is divided into two   parts $\Omega^-_h$ and $\Omega^+_h$  by $\Gamma_h$, which are the approximation of $\Omega^-$ and $\Omega^+$ respectively.

The element in $\mathcal{T}_h$ can be
categorized  into two types:  regular elements and interface elements.  An element $T$ is called interface element if it has exactly two vertices on $\Gamma$; otherwise, it is called regular
element. The set of all interface elements is denoted by $\mathcal{T}_h^*$.   For each element $T \in \mathcal{T}_h^*$, let $T^- = T\cap \Omega^-$ and $T^+=T\cap \Omega^+$.
Since $\Gamma$ is $C^2$, one has
\begin{equation*}
 |T^-| \lesssim h^3_T,  \text{     or        }   |T^+| \lesssim h^3_T,
\end{equation*}
as shown in \cite{ChenZou1998}.

For each edge $e$ on $\Gamma_h$, define a projection $\mathcal{P}_h$  \cite{BrambleKing1996, WeiChenHuangZheng2014} from $e$ to $\Gamma$ as
\begin{equation}
\label{equ:proj}
\mathcal{P}_h(z) = z + d(z)n_h, \quad \forall z \in e,
\end{equation}
where $n_h$ is the unit  normal vector of $e$ pointing from $\Omega^-$ to $\Omega^+$ and $d(z)$ is the sign distance function between $z$ and $\Gamma$ along
$n_h$. Note that $\mathcal{P}_h$ is a point in $\Gamma$  for each $z \in e$.  According to \cite{WeiChenHuangZheng2014, BrambleKing1996}, the projection
  $\mathcal{P}_h$ and its inverse are well defined when the length of $e$ is small enough.

Let $V_h$ be the continuous linear finite element space and $V_{h,0} = V_h \cap H^1_0(\Omega)$. We approximate the diffusion coefficient $\beta$ by $\beta_h$
with $\beta_h|_{T} = \beta^-$ if $T\in \Omega^-_h$ and $\beta_h|_{T} = \beta^+$ if $T\in \Omega^+_h$. Then the linear finite element approximation of the variational problem
\cref{equ:var} is to find $u_h \in V_{h,0} $ such that
\begin{equation}\label{equ:fem}
(\beta_h \nabla u_h, \nabla v_h)  = (f, v_h)  - \langle g_h, v_h\rangle_{\Gamma_h}, \quad \forall v_h \in V_{h,0},
\end{equation}
where $g_h= g(P_h(z))$ and $\langle\cdot, \cdot\rangle$ is $L^2$-inner product of $L^2(\Gamma_h)$. Moreover, \cite{ChenZou1998, Xu1982, Zheng2016}
proved the following convergence results for the finite element approximation \eqref{equ:fem}.

\begin{theorem}
 \label{thm:approx}
 Let $u$ and $u_h$ be the solution to \cref{equ:var} and \cref{equ:fem} respectively, then we have
\begin{align}
 & \|\nabla u- \nabla u_h\|_{0, \Omega}\lesssim h|\log h|^{1/2}(\|f\|_{0, \Omega} + \|g\|_{2, \Gamma}), \label{equ:h1err}\\
  & \|u-u_h\|_{0, \Omega}\lesssim h^2|\log h|^{1/2}(\|f\|_{0, \Omega} + \|g\|_{2, \Gamma}). \label{equ:l2err}
\end{align}
\end{theorem}
%\begin{rem}
 Note that the error estimate \cref{equ:h1err} is nearly optimal due to the existence of $|\log h|^{1/2}$.
 %In the section, we design a post-processing for the interface problem which will produce a nearly superconvergence results.
%\end{rem}

\section{Gradient recovery for elliptic interface problem}
In this section,  we first summarize the polynomial preserving recovery(PPR) method
 proposed by Zhang and Naga in \cite{ ZhangNaga2005, NagaZhang2004, NagaZhang2005} for finite element
 approximation of standard elliptic problem, then based on which, we propose a novel gradient recovery method for elliptic interface problem.
%Moreover, we present one example to illustrate the difference between the proposed method and standard PPR.
\subsection{Polynomial preserving recovery}
%Suppose $A$ is a subset of $\Omega$ and $\mathcal{T}^{A}_h$ is a conforming shape regular triangulation of $A$.
%The set of all mesh vertices  and edges are denoted by  $\mathcal{N}_h^A$ and $\mathcal{E}^A_h$, respectively.  A edge is called
%boundary edge if it is the edge of only one triangle in $\mathcal{T}^{A}_h$. The set of all boundary edges is denoted by $\mathcal{E}^{A,b}_h$.
% A vertex is called boundary vertex if it is one of ending point of a boundary edge $E_h \in \mathcal{E}^{A,b}_h$.
% We want to remark that the boundary edge $E_h \in \mathcal{E}^{A,b}$ or boundary vertice $z \in \mathcal{N}^{A,b}$ may not lie on the real
% boundary of $\Omega$.  For example, if $A = \Omega^-_h$ and $\Gamma \cap \partial\Omega = \emptyset$, boundary edge $E_h \in \mathcal{E}^{A,b}$ ( or boundary vertice $z \in \mathcal{N}^{A,b}$ ) is  just the edge (or vertex)  lie on the approximate interface $\Gamma_h$.

 For any vertex $z$ and $n\in \mathbb{Z}^{+}$,  let $\mathcal{L}(z, n)$ denote the union of  elements
in the first $n$ layers around $z$, i.e.,
\begin{equation}
\mathcal{L}(z, n) := \bigcup \left\{\tau: \tau \in \mathcal{T}_h, \tau \cap \mathcal{L}(z, n-1)\neq \emptyset\right\},
\label{equ:patch}
\end{equation}
where $\mathcal{L}(z, 0):= \{z\}$.

%Let $V_h$ be the continuous linear finite element space defined  on $\mathcal{T}_h$.
The set of all mesh vertices  and edges are denoted by  $\mathcal{N}_h$ and $\mathcal{E}_h$ respectively.
 The  standard Lagrange basis of $V_h$ is denoted by
 $\{\phi_z: z\in \mathcal{N}_h\}$
with $\phi_z(z')=\delta_{zz'}$ for all $z, z'\in \mathcal{N}_h$.
Let us introduce  $G_h: V_h\rightarrow V_h\times V_h$  as the PPR gradient recovery operator.
For any vertex  $z$, let  $\mathcal{K}_z$ be a patch of elements around $z$.
Select all nodes in $\mathcal{N}_h\cap \mathcal{K}_z$ as sampling
points  and fit a polynomial $p_z\in \mathbb{P}_{k+1}(\mathcal{K}_z)$ in  the least
square sense at those sampling points, i.e.
\begin{equation}
  p_z=\arg \min_{p\in \mathbb{P}_{k+1}(\mathcal{K}_z)}
    \sum_{\tilde{z}\in\mathcal{N}_h\cap \mathcal{K}_z}(u_h-p)^2(\tilde{z}).
    \label{least}
\end{equation}
Then the recovered gradient at $z$ is defined as
\begin{equation*}
    (G_hu_h)(z) = \nabla p_z(z).
\end{equation*}

After obtaining recovered gradient value at all nodal points,    we define recovered gradient  $G_h$ on the whole domain by
\begin{equation}
\label{equ:rg}
G_hu_h:= \sum_{z\in \mathcal{N}_h}(G_hu_h)(z)\phi_z.
\end{equation}

\begin{rem}\label{rem:gr}
 If $u_h$ is a function in $V_h$,  then $\nabla u_h$ is a piecewise constant function and hence is discontinuous on $\Omega$.  However, the recovered
 gradient $G_hu_h$ is a continuous piecewise linear function. In that sense, $G_h$ can be viewed as a smoothing operator to smooth a discontinuous piecewise constant
 function into a continuous piecewise linear function.
\end{rem}

To complete the definition of {PPR}, one needs to define $\mathcal{K}_z$. If $z$ is an interior vertex,  $\mathcal{K}_z$
is defined as the smallest $\mathcal{L}(z, n)$  that guarantees the uniqueness of $p_z$ in \eqref{least} \cite{ZhangNaga2005,NagaZhang2004, NagaZhang2005}.
In the case that $z \in \mathcal{N}_h\cup \partial \Omega$, let $n_0$ be the smallest positive integer such that $\mathcal{L}(z, n_0)$
has at least one interior mesh vertex.   Then,  we define
\begin{equation*}
\mathcal{K}_z = \mathcal{L}(z, n_0) \cup \{ \mathcal{K}_{\tilde{z}}: \tilde{z}\in \mathcal{L}(z, n_0) \text{ and }
\tilde{z} \text{ an interior vertex} \}.
\end{equation*}

%\begin{remark}
%     It was proved in \cite{naga2004} that certain rank condition and geometric condition
%    guarantee the uniqueness of $p_z$ in \eqref{least}.
%\end{remark}
\begin{rem}
    In order to avoid numerical instability,  a discrete least squares fitting
    process is carried out on a reference  patch $\omega_{z}$.
\end{rem}

The PPR gradient recovery operator $G_h$ has the following properties, as proved in \cite{ZhangNaga2005, NagaZhang2004, NagaZhang2005}:
\begin{enumerate}
   \item [I.] $G_h$ is a linear operator.
    \item [II.]  $G_h$ preserves quadratic polynomials.   Consequently, $G_h$  enjoys the approximation property
    \begin{equation}\label{equ:pprapprox}
        \|\nabla u - G_h u_I\|\lesssim h^2|u|_{3, \Omega}, \quad \forall u \in H^3(\Omega),
\end{equation}
	where $u_I$ is the linear interpolation of $u$ in $V_h$.
    \item [III.] $\|G_hv_h\|_{0, \Omega}\lesssim \|\nabla v_h\| _{0, \Omega},
	\forall v_h \in V_h$.
\end{enumerate}

\subsection{Immersed polynomial preserving recovery operator}
As we mentioned in \cref{rem:gr}, standard PPR can be viewed a smoothing operator.  However, $\nabla u$ is  discontinuous across the interface $\Gamma$ in elliptic interface problem, and thus standard PPR won't work since it provides continuous gradient approximation. Noticing that,
 although $u$ have low global regularity due to existence of interface, $u|_{\Omega^-}$ (or $u|_{\Omega^+}$) is smooth, which motivates us to recover piecewise continuous gradient approximation instead.

Let $\Omega_h$ be a body-fitted triangulation introduced in Subsection~\ref{ssec:fem}.  The approximate interface $\Gamma_h$ divides the triangulation $\Omega_h$ into two
disjoint sets:
\begin{align}
 &\mathcal{T}^-_h:=\left\{  T\in \mathcal{T}_h| \text{ all three vertices of  } T \text{ are in  }\overline{\Omega^-}  \right\},\label{equ:mmesh}\\
 & \mathcal{T}^+_h:=\left\{  T\in \mathcal{T}_h| \text{ all three vertices of  } T \text{ are in  }\overline{\Omega^+}  \right\} . \label{equ:pmesh}
\end{align}
Let $V_h^-$ and $V_h^-$ be the continuous linear finite element spaces defined on $\mathcal{T}_h^-$ and $\mathcal{T}_h^+$ respectively.

Denote the PPR gradient recovery operator on $V_h^-$ by $G_h^-$.  Then $G_h^-$ is a linear bounded operator from $V_h^-$ to $V_h^-\times V_h^-$.
Similarly, let $G_h^+$ be PPR gradient recovery operator from $V_h^+$ to $V_h^+\times V_h^+$. Then, for any $u_h \in V_h$, we choose the global gradient recovery operator $G_h^I : V_h \rightarrow (V_h^-\cup V_h^+)\times (V_h^-\cup V_h^+)$ as
\begin{equation}
(G_h^I u_h) (z)=
\left\{
\begin{array}{ccc}
   (G_h^- u_h) (z)  &  \text{if } z\in \overline{\Omega^-_h}, \\
   (G_h^+ u_h) (z)  &   \text{if } z\in \overline{\Omega^+_h}.
\end{array}
\right.
\end{equation}
Specifically, we define $(G_h^Iu_h)(z)$ according to the location of $z$:
 \begin{itemize}
\item [Case 1.]  If   $z$ is far from the approximate interface $\Gamma_h$, $(G_h^Iu_h)(z)$ is the standard PPR gradient recovery at $z$.
\item [Case 2.]  If   $z$ is close to the approximate interface $\Gamma_h$, $(G_h^Iu_h)(z)$ is given by fitting a quadratic polynomial using sampling points only from $\mathcal{T}_h^-$ or only from  $\mathcal{T}_h^+$.
\item [Case 3.] If $z$ is on approximate interface $\Gamma_h$, $(G_h^Iu_h)(z)$ is given by two values: one by $(G_h^-u_h)(z)$ and the other by $(G_h^+u_h)(z)$.
\end{itemize}
We call $G_h^I$ as immersed polynomial preserving recovery (IPPR) operator.

\begin{rem}
 Given a function $u_h$ in $V_h$,  $G_h^Iu_h$ is not a function in $V_h\times V_h$, since it is two-valued on the approximate interface
 $\Gamma_h$ as described in Case 3.  %For any vertex $z$ on $\Gamma_h$, $(G_h^Iu_h)(z)$ takes two values: one given by $(G_h^-u_h)(z)$ and the other one si $(G_h^+u_h)(z)$.
\end{rem}

\begin{rem}
 The choice of $G_h^I$ is not limited to PPR gradient recovery operator, and in fact it can be any local gradient recovery operator such as weighted averaging gradient recovery operator \cite{AinsworthOden2000} and SPR gradient recovery operator \cite{ZZ1992a, ZZ1992b}.
One can also use different gradient recovery operators on two subdomains separated by $\Gamma_h$, for example, PPR gradient recovery operator on $V_h^-$ and SPR gradient recovery operator on $V_h^+$. We shall use PPR gradient recovery operator in this paper for convenience.
\end{rem}

\begin{rem}
The IPPR operator can be generalized to the case when the domain $\Omega$ is divided into serval subdomains by defining the gradient recovery operator piecewisely in each subdomain.
\end{rem}

Moreover, we have the following approximation estimate for the IPPR gradient recovery operator $G_h^I$.
\begin{theorem}\label{thm:grbd}
 Let $G_h^I: V_h \rightarrow (V_h^-\cup V_h^+)\times (V_h^-\cup V_h^+)$ be the IPPR gradient recovery operator. Given
 $u\in H^{3}(\Omega^-\cup\Omega^+)$, then
 \begin{equation}
\|G_h^Iu_I - \nabla u \|_{0, \Omega} \lesssim h^{2}\|u\|_{3, \Omega^-\cup\Omega^+},
\label{equ:grbd}
\end{equation}
where $u_I$ is interpolation of $u$ into linear finite element space $V_h$.
\end{theorem}
\begin{proof}
 Notice that $G_h^-$ and $G_h^+$ are the standard PPR gradient recovery operators.   Formula \cref{equ:pprapprox} implies that
\begin{displaymath}
  \|G_h^- u_I- \nabla u \|_{0, \Omega_h^-} \lesssim h^2 \|u\|_{3, \Omega^-}, \text{  and   }  \|G_h^+ u_I- \nabla u \|_{0, \Omega_h^+} \lesssim h^2 \|u\|_{3, \Omega^+}.
\end{displaymath}
Therefore,
\begin{align*}
 \|G_h^Iu_I - \nabla u\|^2_{0, \Omega} &  =   \|G_h^- u_I- \nabla u \|_{0, \Omega_h^-} ^2+   \|G_h^+ u_I- \nabla u \|_{0, \Omega_h^+}^2\\
 &\lesssim h^4 \|u\|_{3, \Omega^-}^2+ h^4 \|u\|_{3, \Omega^+}^2\\
 &\lesssim h^4 \|u\|_{3, \Omega^-\cup\Omega^+}.
\end{align*}
\end{proof}

\section{Superconvergence Analysis}
In this section, we prove that the IPPR gradient recovery method has superconvergence for both mildly unstructured mesh and adaptive mesh.

%Using an supercloseness result which can proved similarly as in \cite{WeiChenHuangZheng2014}, we firstly prove the superconvergence results
%on mildly unstructured mesh.  Then we establish the supercloseness result on adaptive refined mesh.  Based on the supercloseness result,
%we can show that our gradient recovery operator is also enjoy superconvergence on adaptive refined mesh and hence it provide a
%asymptotically exact a posteriori error estimator for interface.

\subsection{Superconvergence on mildly unstructured mesh} We first introduce a definition on mesh structure.
\begin{definition}
1. Two adjacent triangles are called to  form an $\mathcal{O}(h^{1+\alpha})$  approximate parallelogram if the lengths of any two opposite edges differ only by $\mathcal{O}(h^{1+\alpha})$.

2. The triangulation $\mathcal{T}_h$ is called to satisfy Condition
$(\sigma,\alpha)$ if
there exist a partition $\mathcal{T}_{h,1} \cup \mathcal{T}_{h,2}$ of $\mathcal{T}_h$
and positive constants $\alpha$ and $\sigma$ such that every
two adjacent triangles in $\mathcal{T}_{h,1}$ form an $\mathcal{O}(h^{1+\alpha})$ parallelogram and
$$
\sum_{T\in {\mathcal{T}_{h,2}}} |T| = \mathcal{O}(h^\sigma).
$$
\end{definition}

Under the above mesh condition, the following supercloseness result holds:
\begin{theorem}
\label{thm:superclose}
Let $u$ be the solution to variational problem \cref{equ:var} and $u_h$ be the finite element solution to \cref{equ:fem}.
 If the body-fitted mesh satisfies Condition $(\sigma,\alpha)$ and $u \in H^1(\Omega) \cap H^3(\Omega^-\cup\Omega^+)
 \cap W^{2, \infty}(\Omega^-\cup\Omega^+)$, then  for any $v_h \in V_{h,0}$,
 \begin{equation}\label{equ:approximation}
 \begin{split}
(\beta_h \nabla(u-u_I), \nabla v_h)  \lesssim& h^{1+\rho}(\|u\|_{3, \Omega^-\cup\Omega^+}
+ \|u\|_{2, \infty,  \Omega^-\cup\Omega^+})|v_h|_{1, \Omega}  \\
&+h^{\frac{3}{2}}\|u\|_{2, \infty,  \Omega^-\cup\Omega^+}|v_h|_{1, \Omega},
\end{split}
\end{equation}
and
  \begin{equation}\label{equ:superclose}
 \begin{split}
 \|\nabla(u_I-u_h)\|_{0, \Omega} \lesssim& h^{1+\rho}(\|u\|_{3, \Omega^-\cup\Omega^+}
+ \|u\|_{2, \infty,  \Omega^-\cup\Omega^+})  \\
&+h^{\frac{3}{2}}(\|u\|_{2, \infty,  \Omega^-\cup\Omega^+}+\|g\|_{0, \infty, \Gamma}),
\end{split}
\end{equation}
where $\rho = \min(\alpha, \frac{\sigma}{2}, \frac{1}{2})$ and $u_I \in V_h$ is the interpolation of $u$.
\end{theorem}
\begin{proof}
The proof is similar to the proof of Theorem 3.6 in \cite{WeiChenHuangZheng2014} where one uses the estimates in \cite{XuZhang2004} instead of
\cite{ChenXu2007}.
\end{proof}

{
\begin{rem}\label{rem:improved}
 If $\Gamma=\Gamma_h$ or the flux jump $g$ given by \eqref{equ:fluxjump} vanishes with $\Gamma$ smooth enough so that $ |T^-|  \lesssim h^4_T $  ( or  $ |T^+|  \lesssim h^4_T $ ),
  then one has the following $O(h^{1+\rho})$ superconvergence results
   \begin{equation}\label{equ:improvedap}
 \begin{split}
(\beta_h \nabla(u-u_I), \nabla v_h)  \lesssim& h^{1+\rho}(\|u\|_{3, \Omega^-\cup\Omega^+}
+ \|u\|_{2, \infty,  \Omega^-\cup\Omega^+})|v_h|_{1, \Omega},
\end{split}
\end{equation}
and
  \begin{equation}\label{equ:improvedsuperclose}
 \begin{split}
 \|\nabla(u_I-u_h)\|_{0, \Omega} \lesssim& h^{1+\rho}(\|u\|_{3, \Omega^-\cup\Omega^+}
+ \|u\|_{2, \infty,  \Omega^-\cup\Omega^+}).
\end{split}
\end{equation}
\end{rem}
}

\begin{rem}
 The difference between \cref{thm:superclose} and Theorem 3.6 in \cite{WeiChenHuangZheng2014} lies in the condition on mesh structure.  The $\mathcal{O}(h^{2\sigma})$
 irregular mesh in \cite{WeiChenHuangZheng2014} is a special case of Condition $(\sigma,\alpha)$.
\end{rem}

\begin{rem}
 For body-fitted mesh generated by Delaunay algorithm,
   the assumption that two adjacent triangles form $\mathcal{O}(h^{1+\alpha})$  approximate parallelogram is violated near boundary and interface, however, the summation of
 area of such triangles is bounded by $ \mathcal{O}(h^\sigma)$, and thus it Condition $(\sigma,\alpha)$.
\end{rem}

The supercloseness result in \cref{thm:superclose} implies the following superconvergent result.
\begin{theorem}\label{thm:superconvergence}
Under the same hypothesis as in \Cref{thm:superclose}, then we have
 \begin{equation}\label{equ:super}
 \begin{split}
 \|\nabla u-G_h^Iu_h \|_{0, \Omega} \lesssim& h^{1+\rho}(\|u\|_{3, \Omega^-\cup\Omega^+}
+ \|u\|_{2, \infty,  \Omega^-\cup\Omega^+}) + \\
&h^{\frac{3}{2}}(\|u\|_{2, \infty,  \Omega^-\cup\Omega^+}+\|g\|_{0, \infty, \Gamma}),
\end{split}
\end{equation}
where $\rho = \min(\alpha, \frac{\sigma}{2}, \frac{1}{2})$.
\end{theorem}
\begin{proof}
 We decompose $\nabla u-G_h^Iu_h $ as  $(\nabla u-G_h^Iu_I)- (G_h^Iu_I - G_h^I u_h ) $, and then by triangle inequality,
 \begin{equation}
\begin{split}
&\| \nabla u-G_h^Iu_h \|_{0, \Omega}\\
 \le &\|\nabla u-G_h^Iu_I\|_{0,\Omega} + \| G_h^Iu_I - G_h^I u_h  \|_{0,\Omega}\\
\le &\|\nabla u- G_h^Iu_I \|_{0,\Omega}+\|G_h^-u_I- G_h^- u_h \|_{0,\Omega_h^-} +\|G_h^+u_I - G_h^+ u_h \|_{0,\Omega_h^+}\\
\le & \|\nabla u- G_h^Iu_I \|_{0,\Omega}+\|\nabla(u_I -  u_h ))\|_{0,\Omega_h^-} +\|\nabla(u_I -  u_h)\|_{0,\Omega_h^+}\\
\le &\|\nabla u- G_h^Iu_I \|_{0,\Omega}+\|\nabla(u_I -  u_h )\|_{0,\Omega} ,\\
\end{split}
\end{equation}
where we have used boundedness of PPR gradient recovery operator $G_h^-$ and $G_h^+$.   The first term can be bounded by
\Cref{thm:grbd} and the second term is estimated in  \Cref{thm:superclose}, which completes the proof.
\end{proof}

{
\begin{rem}\label{rem:improvedsuper}
 Under the same assumptions in \cref{rem:improved}, one has the following  improved  superconvergence results
  \begin{equation}\label{equ:improvedsuper}
 \begin{split}
 \|\nabla u-G_h^Iu_h\|_{0, \Omega} \lesssim& h^{1+\rho}(\|u\|_{3, \Omega^-\cup\Omega^+}
+ \|u\|_{2, \infty,  \Omega^-\cup\Omega^+}).
\end{split}
\end{equation}
\end{rem}
}

\subsection{Superconvergence on adaptive mesh}  In this subsection, for simplicity, we assume that the interface $\Gamma$ does not cut through  any element $T\in \mathcal{T}_h$, i.e. $\Gamma=\Gamma_h$, which implies $\Omega^- = \Omega^-_h$  and $\Omega^+ = \Omega^+_h$.
Furthermore, we assume that the solution $u$ to \cref{equ:var} has a single singularity on the interface
$\Gamma$, and without loss of generality, we assume that the singularity is at the origin and $u^{-}$ (or $u^+$) can be decomposed into a smooth
part $w^-$ (or $w^+$)  and singular part  $v^-$ (or $v^+$), i.e.,
\begin{equation}\label{equ:decomp}
u^- = w^- + v^-, \text{   and    }   u^+ = w^+ + v^+,
\end{equation}
where
\begin{equation}
 \left|\frac{\partial^m w^{-}}{\partial x^i\partial y^{m-i}}\right| \lesssim 1,  \text{    and      }
   \left|\frac{\partial^m v^{-}}{\partial x^i\partial y^{m-i}}\right| \lesssim r^{\delta -m},
\end{equation}
and
\begin{equation}
 \left|\frac{\partial^m w^{+}}{\partial x^i\partial y^{m-i}}\right| \lesssim 1,  \text{    and      }
   \left|\frac{\partial^m v^{+}}{\partial x^i\partial y^{m-i}}\right| \lesssim r^{\delta -m},
\end{equation}
with $r = \sqrt{x^2+y^2}$ and $0< \delta <2$ being a constant.

%Remark that these assumptions are true for Kellogg problem \cite{Kellogg1974} which is a typical model for the study of adaptive finite element method on interface problem \cite{ChenDai2002,CaiZhang2009, MorinNochetto2003}.

%reasonable for elliptic interface problem.
%For instance, such hypothesis is true for Kellogg problem \cite{Kellogg1974} which is a model problem for studying adaptive finite element method for
% interface problem \cite{ChenDai2002,CaiZhang2009, MorinNochetto2003}.

For any edge $e$ of the mesh $\mathcal{T}_h$, let $h_e$ be the length of the edge and $r_e$ be the distance from the origin to the midpoint of $e$.
  If $e$ is an interior edge, denote $\Omega_e$ to be
the patch of $e$ consisting of two triangles sharing the edge $e$. In addition,  let  $\underline{h} \simeq \min_{T\in \mathcal{T}_h} h_T$
and $N$ be the number of vertices of $\mathcal{T}_h$.  To get superconvergence, one also needs the following restriction on mesh structure.
\begin{definition}
 The triangulation $\mathcal{T}_h$ is called to satisfy \text{Condition} $(\alpha, \sigma, \mu)$ if  there exist constants $\alpha>0$, $\sigma \ge 0$,
 and $\mu >0$ such that the interior edge can separated into two parts $\mathcal{E}_h = \mathcal{E}_{1,h}\oplus \mathcal{E}_{2,h}$:
 $\Omega_e$ forms an $\mathcal{O}\left(h^{1+\alpha}_e/r_e^{\alpha+\mu(1-\alpha)}\right)$ parallelogram for $\forall e \in \mathcal{E}_{1,h}$
 and the number of edges in $\mathcal{E}_{2,h}$ satisfies $\#\mathcal{E}_{2,h}\lesssim N^{\sigma}$.
\end{definition}

Note that the above mesh condition is a practical assumption for adaptive mesh as shown in \cite{WuZhang2007}.  In addition, we assume $h_T \simeq r_T^{1-\mu}\underline{h}^{\mu}$
for any $T\in \mathcal{T}_h$, and then we can establish the following supercloseness result on adaptive mesh.
 %For adaptive refined mesh interface-fitted mesh, the mesh $\mathcal{T}_h$ satisfies the following mesh condition \cite{WuZhang2007}:
\begin{theorem}\label{thm:supercloseafem}
Let $u$ be the solution to variational problem \cref{equ:var} and $u_h$ be the finite element solution to \cref{equ:fem}.
 Suppose  adaptive refined  mesh $\mathcal{T}_h$ satisfies \text{Condition} $(\alpha, \sigma, \delta/2)$, and $h_T\simeq r_T^{1-\delta/2}
 \underline{h}^{\delta/2}$ for any $T \in \mathcal{T}_h$, then for any $v_h \in V_{h,0}$,
 \begin{equation}\label{equ:afemapprox}
(\beta_h \nabla(u-u_I), \nabla v_h)  \lesssim \frac{1+(\ln N) ^{1/2}}{N^{1/2+\rho}}\|\nabla v_h\|_{0, \Omega},
\end{equation}
and
  \begin{equation}\label{equ:afemsuperclose}
 \begin{split}
 \|\beta^{1/2}(\nabla u_I-\nabla u_h)\|_{0, \Omega} \lesssim& \frac{1+(\ln N) ^{1/2}}{N^{1/2+\rho}}\|\nabla v_h\|_{0, \Omega},
\end{split}
\end{equation}
where $\rho = \min\left(\frac{\alpha}{2}, \frac{1-\sigma}{2}\right)$ and $u_I$ is the linear interpolation of $u$.
\end{theorem}
\begin{proof}
First, one can decompose $(\beta_h \nabla(u-u_I), \nabla v_h)$  as
 \begin{equation*}
\begin{split}
   &(\beta_h \nabla(u-u_I), \nabla v_h) \\
%    =& (\beta \nabla(u-u_I), \nabla v_h) \\
  = &(\beta^-( \nabla u^--u_I^-), \nabla v_h)_{\Omega^-} + (\beta^+( \nabla u^+-u_I^+), \nabla v_h)_{\Omega^+}\\
%  \lesssim &\frac{1+(\ln N) ^{1/2}}{N^{1/2+\rho}}\|\nabla v_h\|_{0, \Omega},\\
  =&I_1+I_2.
    \end{split}
\end{equation*}
Lemma 3.3 in \cite{WuZhang2007} implies the following estimates for $I_1$ and $I_2$,
\begin{equation*}
\begin{split}
 I_1 \lesssim  &\frac{1+(\ln N) ^{1/2}}{N^{1/2+\rho}}\|\nabla v_h\|_{0, \Omega},\\
 I_2 \lesssim  &\frac{1+(\ln N) ^{1/2}}{N^{1/2+\rho}}\|\nabla v_h\|_{0, \Omega},
\end{split}
\end{equation*}
which completes the proof of \cref{equ:afemapprox}.  By \cref{equ:var,equ:fem} and noticing that
$\Gamma=\Gamma_h$ and $\beta = \beta_h$, we have
           \begin{equation*}
(\beta_h \nabla(u_I-u_h), \nabla v_h)  = (\beta_h \nabla(u_I-u), \nabla v_h).
\end{equation*}
Taking $v_h = u_I-u_h$ gives \cref{equ:afemsuperclose}.
\end{proof}

Before presenting our main superconvergent theorem on adaptive refined mesh, we need to estimate gradient recovery operator analogous to \Cref{thm:grbd}.
\begin{theorem}\label{thm:approxafem}
 Assume that $h_T \simeq r_T^{1-\delta/2}\underline{h}^{\delta/2}$ for any $T \in \mathcal{T}_h$, then
 \begin{equation}\label{equ:approxafem}
\|\nabla u - G_h^Iu_I \|_{0, \Omega} \lesssim \frac{1+(\ln N)^{1/2}}{N}.
\end{equation}
\end{theorem}
\begin{proof}
The definition of $G_h^I$ produces
 \begin{equation}\label{equ:grdecomp}
\|\nabla u - G_h^Iu_I\|_{0, \Omega}^2  =  \|\nabla u^- G_h^-u_I^- \|_{0, \Omega^-}^2 +   \|\nabla u^+ - G_h^+u_I^+ \|_{0, \Omega^+}^2.
\end{equation}
Note that $u^-$ (or $u^+$) has the decomposition \cref{equ:decomp} and $G_h^-$ (or  $G_h^+$ ) is the standard PPR gradient recovery operator.
Lemma 5.2 in \cite{WuZhang2007} implies
\begin{align}
  \|\nabla u^- - G_h^-u_I^-\|_{0, \Omega^-} \lesssim  & \frac{1+(\ln N)^{1/2}}{N},\label{equ:mbdafem}\\
  \|\nabla u^+ - G_h^+u_I^+ \|_{0, \Omega^+} \lesssim  &\frac{1+(\ln N)^{1/2}}{N}\label{equ:pbdafem}.
\end{align}
Combing \cref{equ:grdecomp,equ:mbdafem,equ:pbdafem} gives \cref{equ:approxafem}.
\end{proof}

Then we can prove the superconvergence as follows.
\begin{theorem}\label{thm:superafem}
 Let $u$ be the solution to variational problem \cref{equ:var} and $u_h$ be the finite element solution to \cref{equ:fem}.
 If the adaptive refined  mesh $\mathcal{T}_h$ satisfies \text{Condition} $(\alpha, \sigma, \delta/2)$ and $h_T\simeq r_T^{1-\delta/2}
 \underline{h}^{\delta/2}$ for any $T \in \mathcal{T}_h$, then
\begin{equation}\label{equ:superafem}
\|\beta^{1/2}(\nabla u- G_h^Iu_h)\|_{0, \Omega} \lesssim  \frac{1+(\ln N) ^{1/2}}{N^{1/2+\rho}},
\end{equation}
where $\rho = \min\left(\frac{\alpha}{2}, \frac{1-\sigma}{2}\right)$.
\end{theorem}
\begin{proof}
 The proof is essentially the same as in \Cref{thm:superconvergence} where one should use \Cref{thm:approxafem,thm:supercloseafem} instead of  \Cref{thm:grbd,thm:superclose}.
\end{proof}

Using the IPPR gradient recovery operator $G_h^I$, one can define a local {\it a posteriori} error
estimator on element $T\in \mathcal{T}_h$ as
\begin{equation}\label{equ:localind}
\eta_T = \|\beta^{1/2}(G_h^Iu_h - \nabla u_h)\|_{0, T},
\end{equation}
and the corresponding global error estimator  as
\begin{equation}\label{equ:globalind}
\eta_h = \left( \sum_{T\in \mathcal{T}_h}\eta_T^2\right)^{1/2}.
\end{equation}

\Cref{thm:superafem} implies that $\eta_h$ is asymptotically exact {\it a posteriori} error estimator for elliptic
 interface problem.
\begin{theorem}\label{thm:asy}
 Let $u$ be the solution to variational problem \cref{equ:var} and $u_h$ be the finite element solution to \cref{equ:fem}.
 Suppose adaptive refined  mesh $\mathcal{T}_h$ satisfies \text{Condition} $(\alpha, \sigma, \delta/2)$ and $h_T\simeq r_T^{1-\delta/2}
 \underline{h}^{\delta/2}$ for any $T \in \mathcal{T}_h$, and if
 \begin{equation}\label{equ:lowerbnd}
\frac{1}{N^{1/2}} \lesssim \|\beta^{1/2}(\nabla(u-u_h))\|_{0, \Omega},
\end{equation}
then
\begin{equation}
 \left| \frac{\eta_h}{ \|\beta^{1/2}(\nabla u-\nabla u_h)\|_{0, \Omega}} - 1 \right|   \lesssim  \frac{1+(\ln N) ^{1/2}}{N^{\rho}},
\end{equation}
where $\rho = \min\left(\frac{\alpha}{2}, \frac{1-\sigma}{2}\right)$.
\end{theorem}

\begin{rem}
 Assumption \cref{equ:lowerbnd} is reasonable according to lower bounded estimates of approximation by finite element spaces
 \cite{Krizek2011, LinXieXu2014}.
\end{rem}

\section{Numerical Examples}
In this section,  we present serval numerical examples to illustrate  the superconvergence
of the IPPR gradient recovery method and confirm the theoretical results given in the previous section.
We also make numerical comparison to standard PPR method \cite{NagaZhang2004, NagaZhang2005, ZhangNaga2005} to show the effectiveness.
%All numerical computations are done by MATLAB.
For convenience, we shall use the following error measurements in all examples:
\begin{eqnarray}\label{eq:notations-1}
De&:=&\|u-u_h\|_{1,\Omega},\quad\quad\quad D^ie:=\|\nabla u_I- \nabla u_h\|_{0, \Omega},\\
D^re&:=&\|\nabla u-G^I_hu_h\|_{0, \Omega},\quad
D^pe:=\|\nabla u-G_hu_h\|_{ 0, \Omega}. \label{eq:notations-2}
 \end{eqnarray}

Remark that all convergence rate will be computed in the degree of freedom (Dof), and
since  $\text{Dof}  \approx h^{-2}$ for a two-dimensional quasi-uniform mesh,  the
corresponding  convergent rate in mesh size $h$ is twice as much as what we present in the tables.

{\bf Example 5.1.} In this example, we  consider  the elliptic interface problem  \cref{equ:model} with a circular interface of radius $r_0 = 0.5$ as studied in \cite{LiLinWu2003}.
 The exact solution is
\begin{equation*}
u(z) =
\left\{
\begin{array}{ll}
    \frac{r^3}{\beta^-}   &  \text{if }   z\in \Omega_-, \\
      \frac{r^3}{\beta^+} + \left( \frac{1}{\beta^-}-\frac{1}{\beta^+} \right)r_0^3&  \text{if } z \in \Omega^+,\\
   \end{array}
\right.
\end{equation*}
where $r = \sqrt{x^2+y^2}$.

\begin{figure}[ht]
\centering
\subfigure[]{%
     \includegraphics[width=0.4\textwidth]{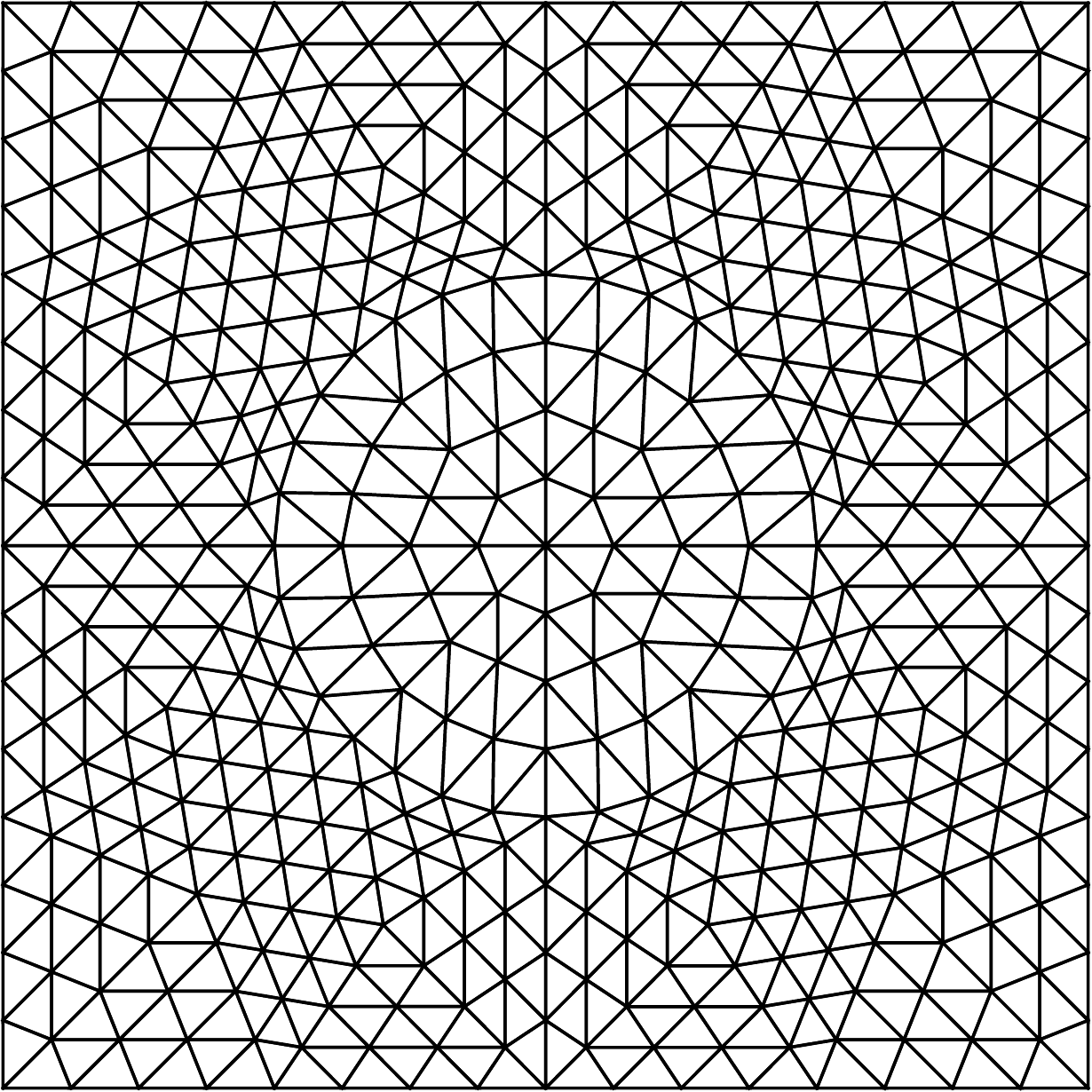}
  \label{fig:circlemesh}}
\quad
\subfigure[]{%
     \includegraphics[width=0.53\textwidth]{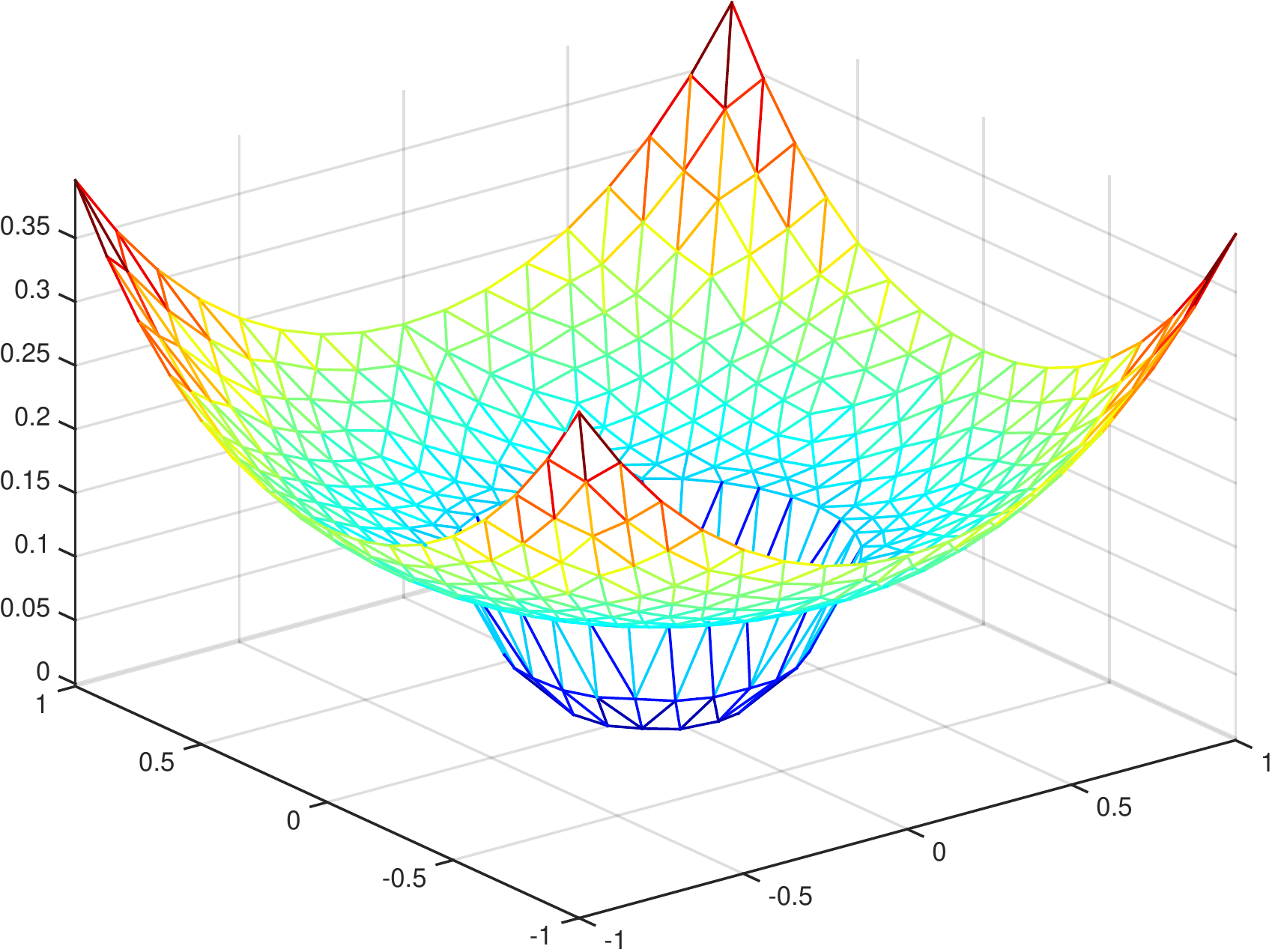}
  \label{fig:circlesol}}
\caption{Finite element mesh and solution for Example 5.1 with $\beta^+=10, \beta^-=1$:  (a)  Body-fitted mesh on the second level; (b) Finite element solution on the body-fitted mesh.}
\label{fig:circle}
\end{figure}

\begin{figure}[ht]
\centering
\subfigure[]{%
     \includegraphics[width=0.45\textwidth]{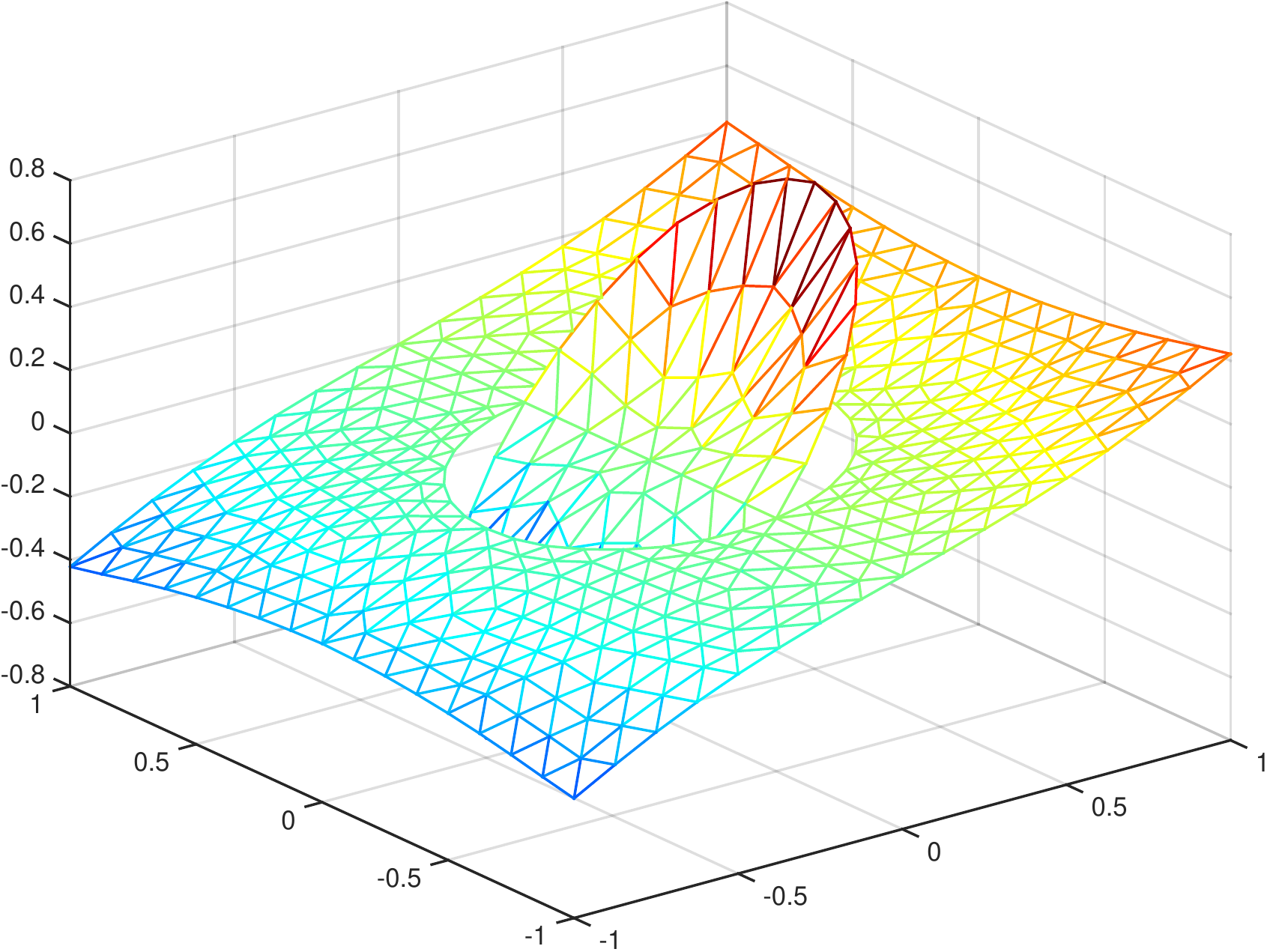}
  \label{fig:circlecx}}
\quad
\subfigure[]{%
     \includegraphics[width=0.45\textwidth]{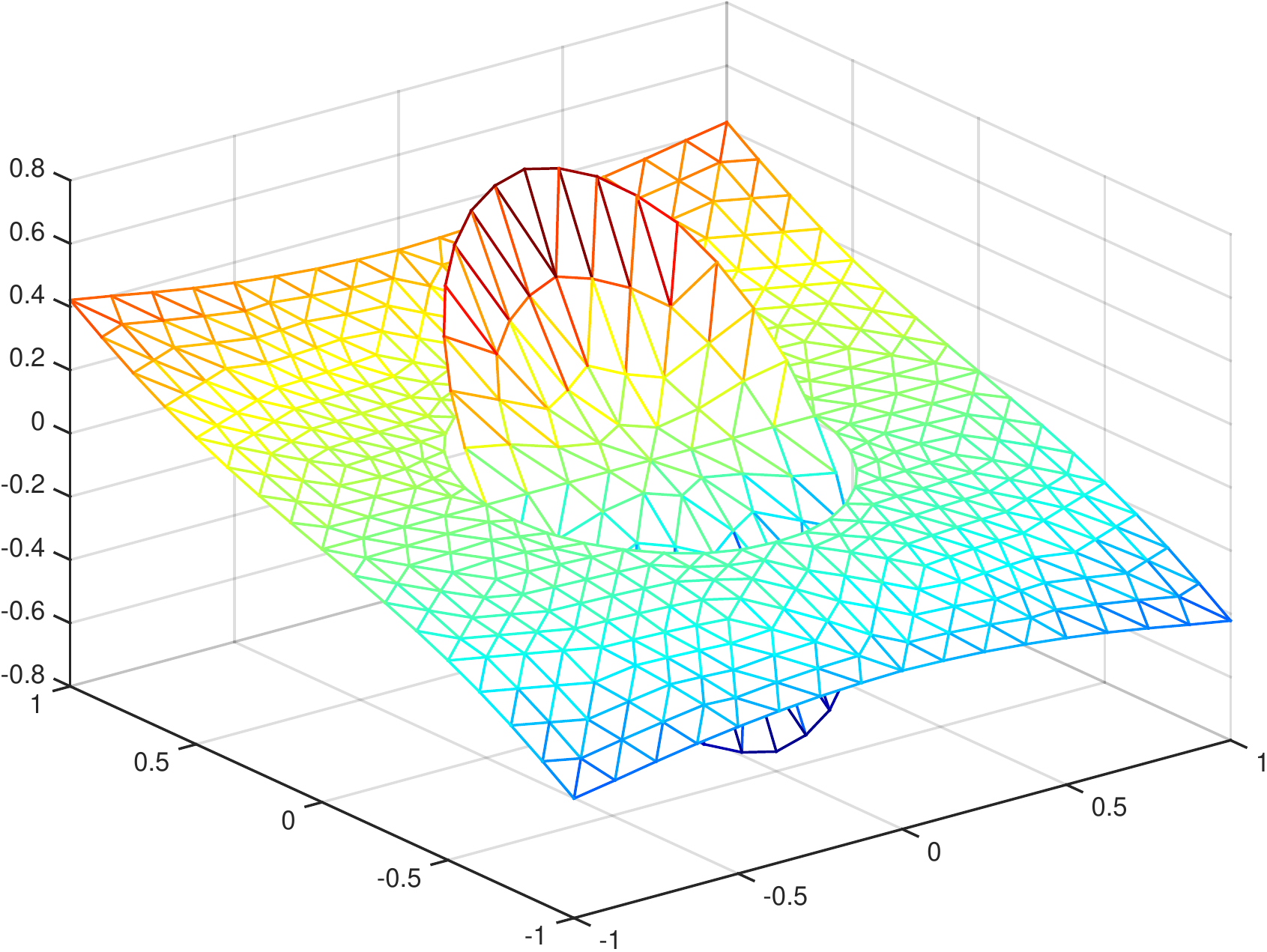}
  \label{fig:circlery}}
\caption{Plot of  recovered gradient for Example 5.1 with $\beta^+=10, \beta^-=1$: (a) $x$-component; (b) $y$-component.}
\label{fig:circlegradient}
\end{figure}

Here we use five different levels of body-fitted meshes generated by  Delaunay mesh generator.
\cref{fig:circlemesh} plots the second level body-fitted mesh and \cref{fig:circlesol} plots the finite element solution on such
a mesh.  \cref{fig:circlegradient} shows the recovered gradient.

\begin{table}[htb!]\label{tab:ex51a}
\centering
\caption{Numerical results for Example 5.1 with $\beta^+=10, \beta^-=1$. }
\begin{tabular}{|c|c|c|c|c|c|c|c|c|c|}
\hline
 Dof & $De$ & order& $D^{i}e$ & order& $D^{r}e$ & order& $D^{p}e$ &  order\\ \hline\hline
 129 &1.35e-01&--&1.63e-02&--&1.34e-01&--&2.44e-01&--\\ \hline
 481 &7.25e-02&0.48&5.00e-03&0.90&2.30e-02&1.34&1.80e-01&0.23\\ \hline
 1857 &3.69e-02&0.50&1.40e-03&0.94&6.82e-03&0.90&1.30e-01&0.24\\ \hline
 7297 &1.85e-02&0.50&3.76e-04&0.96&1.84e-03&0.96&9.24e-02&0.25\\ \hline
 28929 &9.29e-03&0.50&9.81e-05&0.97&4.78e-04&0.98&6.52e-02&0.25\\ \hline
\end{tabular}
\end{table}

\begin{table}[htb!]\label{tab:ex51b}
\centering
\caption{Numerical results for Example 5.1 with $\beta^+=1000, \beta^-=1$. }
\begin{tabular}{|c|c|c|c|c|c|c|c|c|c|}
\hline
 Dof & $De$ & order& $D^{i}e$ & order& $D^{r}e$ & order& $D^{p}e$ &  order\\ \hline\hline
 129 &1.25e-01&--&1.51e-02&--&1.47e-01&--&2.76e-01&--\\ \hline
 481 &6.76e-02&0.47&4.63e-03&0.90&2.29e-02&1.41&2.01e-01&0.24\\ \hline
 1857 &3.45e-02&0.50&1.29e-03&0.95&6.78e-03&0.90&1.45e-01&0.24\\ \hline
 7297 &1.73e-02&0.50&3.39e-04&0.98&1.83e-03&0.96&1.03e-01&0.25\\ \hline
 28929 &8.68e-03&0.50&8.70e-05&0.99&4.75e-04&0.98&7.24e-02&0.25\\ \hline
\end{tabular}
\end{table}

\begin{table}[htb!]\label{tab:ex51c}
\centering
\caption{Numerical results for Example 5.1 with $\beta^+=1000000, \beta^-=1$. }
\begin{tabular}{|c|c|c|c|c|c|c|c|c|c|}
\hline
 Dof & $De$ & order& $D^{i}e$ & order& $D^{r}e$ & order& $D^{p}e$ &  order\\ \hline\hline
 129 &1.25e-01&--&1.51e-02&--&1.47e-01&0.00&2.76e-01&--\\ \hline
 481 &6.76e-02&0.47&4.63e-03&0.90&2.29e-02&1.41&2.01e-01&0.24\\ \hline
 1857 &3.45e-02&0.50&1.29e-03&0.95&6.78e-03&0.90&1.45e-01&0.24\\ \hline
 7297 &1.73e-02&0.50&3.39e-04&0.98&1.83e-03&0.96&1.03e-01&0.25\\ \hline
 28929 &8.68e-03&0.50&8.70e-05&0.99&4.75e-04&0.98&7.25e-02&0.25\\ \hline
\end{tabular}
\end{table}

\begin{table}[htb!]\label{tab:ex51d}
\centering
\caption{Numerical results for Example 5.1 with $\beta^+=1, \beta^-=1000000$. }
\begin{tabular}{|c|c|c|c|c|c|c|c|c|c|}
\hline
 Dof & $De$ & order& $D^{i}e$ & order& $D^{r}e$ & order& $D^{p}e$ &  order\\ \hline\hline
 129 &5.27e-01&--&6.02e-02&--&1.84e-01&0.00&3.51e-01&--\\ \hline
 481 &2.64e-01&0.53&1.70e-02&0.96&3.06e-02&1.36&2.14e-01&0.38\\ \hline
 1857 &1.32e-01&0.51&4.66e-03&0.96&8.06e-03&0.99&1.48e-01&0.28\\ \hline
 7297 &6.60e-02&0.51&1.26e-03&0.96&2.10e-03&0.98&1.03e-01&0.26\\ \hline
 28929 &3.30e-02&0.50&3.37e-04&0.96&5.40e-04&0.98&7.25e-02&0.26\\ \hline
\end{tabular}
\end{table}

%Using notations in \cref{eq:notations-1,eq:notations-2},
\cref{tab:ex51a,tab:ex51b,tab:ex51c,tab:ex51d} show the numerical results for four typical different
jump ratios: $\beta^-/\beta^+ = 1/10$ (moderate jump), $\beta^-/\beta^+ = 1/1000$ (large jump), $\beta^-/\beta^+ = 1/1000000$ (huge jump),
and $\beta^-/\beta^+ = 1000000$ (huge jump).   In all different cases,  optimal $\mathcal{O}(h)$ convergence can be observed for $H^1$-semi error
of finite element solution as given in \cref{thm:approx}.
 Notice that in this example, $g=0$ and $\Gamma$ is arbitrary smooth.
As discussed in \cref{rem:improved}, one can have the supercloseness of $\mathcal{O}(h^2)$  as observed in Column 5 of \cref{tab:ex51a,tab:ex51b,tab:ex51c,tab:ex51d}.
One also notices that, for the convergence rate of gradients, IPPR ($D^re$) superconverges at the order of $\mathcal{O}(h^2)$ while PPR ($D^pe$) converges suboptimally at the order of $\mathcal{O}(h^{0.5})$.

{\bf Example 5.2.}  In this example,  we consider the flower-shape interface problem as studied in \cite{MuWang2013, ZhouWei2006}.
The computational domain is $(-1, 1)\times (-1,1)$.
The interface curve $\Gamma$ in polar coordinate is given by
\begin{equation*}
 r = \frac{1}{2} + \frac{\sin(5\theta)}{7},
\end{equation*}
which contains both convex and concave parts.
The diffusion coefficient is piecewise constant with $\beta^-=1$ and $\beta^+=10$.   The  right hand function $f$ in \eqref{equ:model} is chosen to match the exact solution
\begin{equation*}
u(z) =
\left\{
\begin{array}{ll}
    e^{(x^2+y^2)}, & \text{if } z \in \Omega^+\\
    0.1(x^2+y^2)^2-0.01\ln(2\sqrt{x^2+y^2}),&  \text{if } z \in \Omega^+,\\
   \end{array}
\right.
\end{equation*}
and the jump conditions at interface \eqref{equ:valuejump}-\eqref{equ:fluxjump} are also provided by the exact solution.

\begin{figure}[ht]
\centering
\subfigure[]{%
     \includegraphics[width=0.4\textwidth]{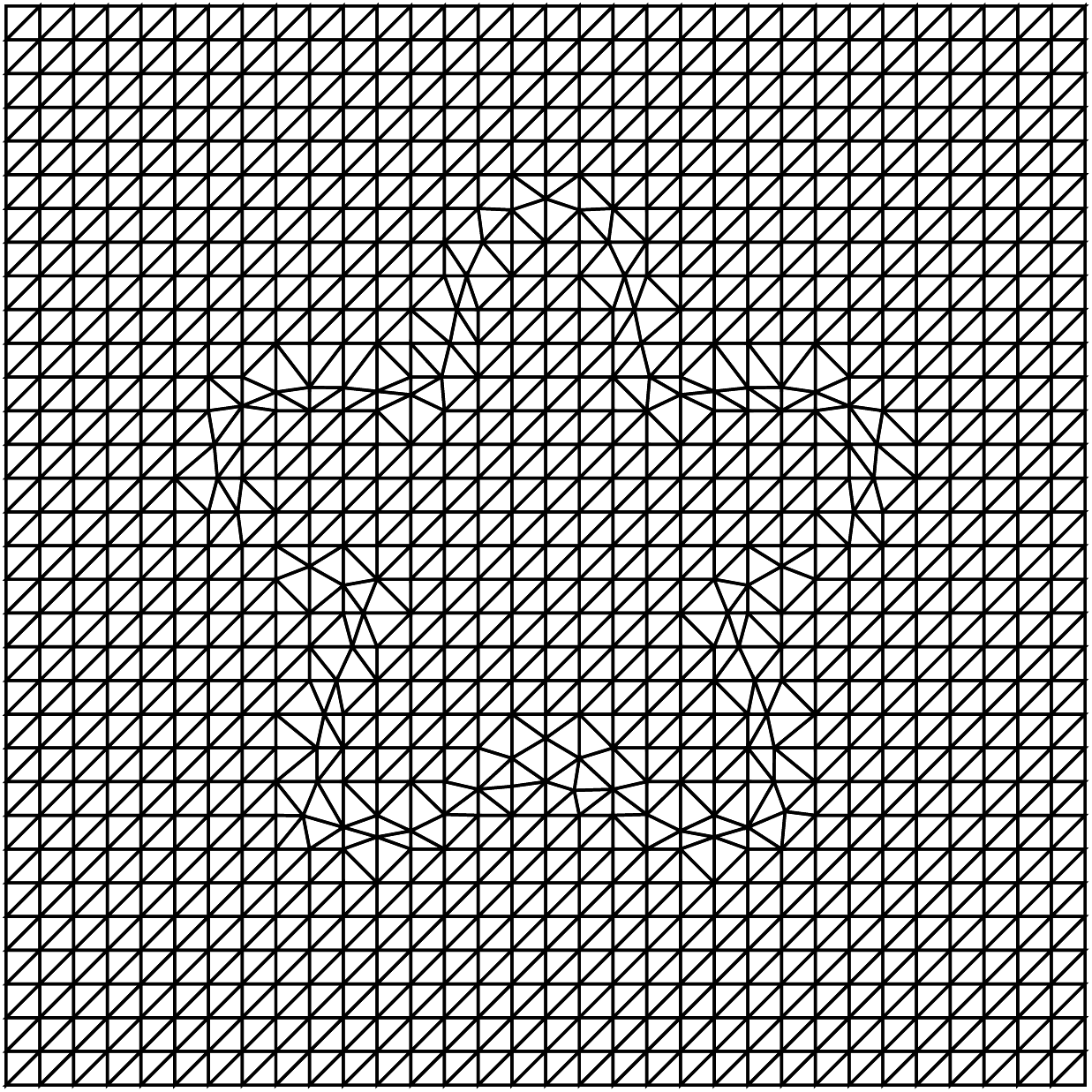}
  \label{fig:flowermesh}}
\quad
\subfigure[]{%
     \includegraphics[width=0.53\textwidth]{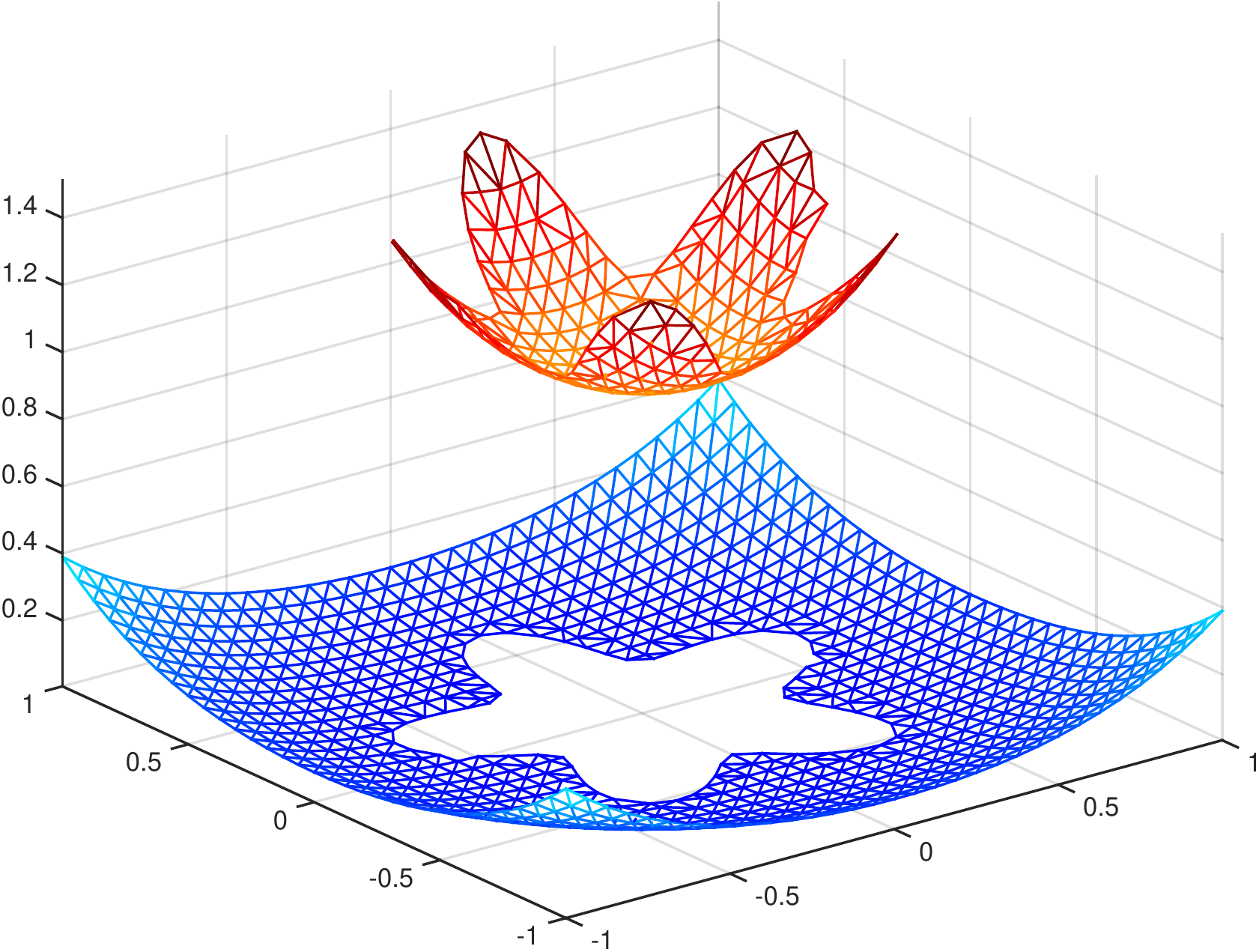}
  \label{fig:flowersol}}
\caption{Finite element mesh and solution for Example 5.2:  (a) Body-fitted mesh on the first level; (b) Finite element solution on the body-fitted mesh.}
\label{fig:flower}
\end{figure}

\begin{figure}[ht]
\centering
\subfigure[]{%
     \includegraphics[width=0.45\textwidth]{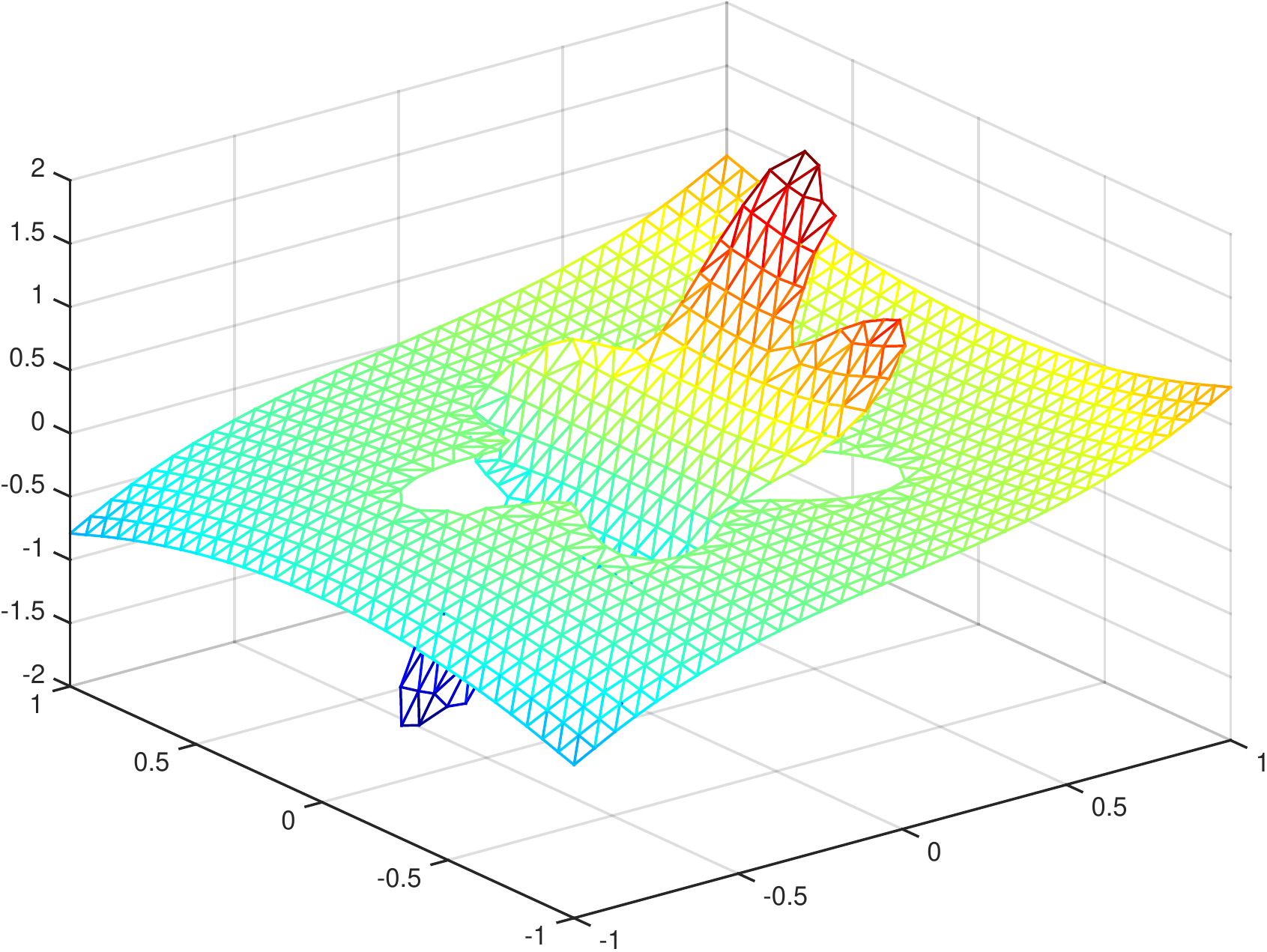}
  \label{fig:flowerrx}}
\quad
\subfigure[]{%
     \includegraphics[width=0.45\textwidth]{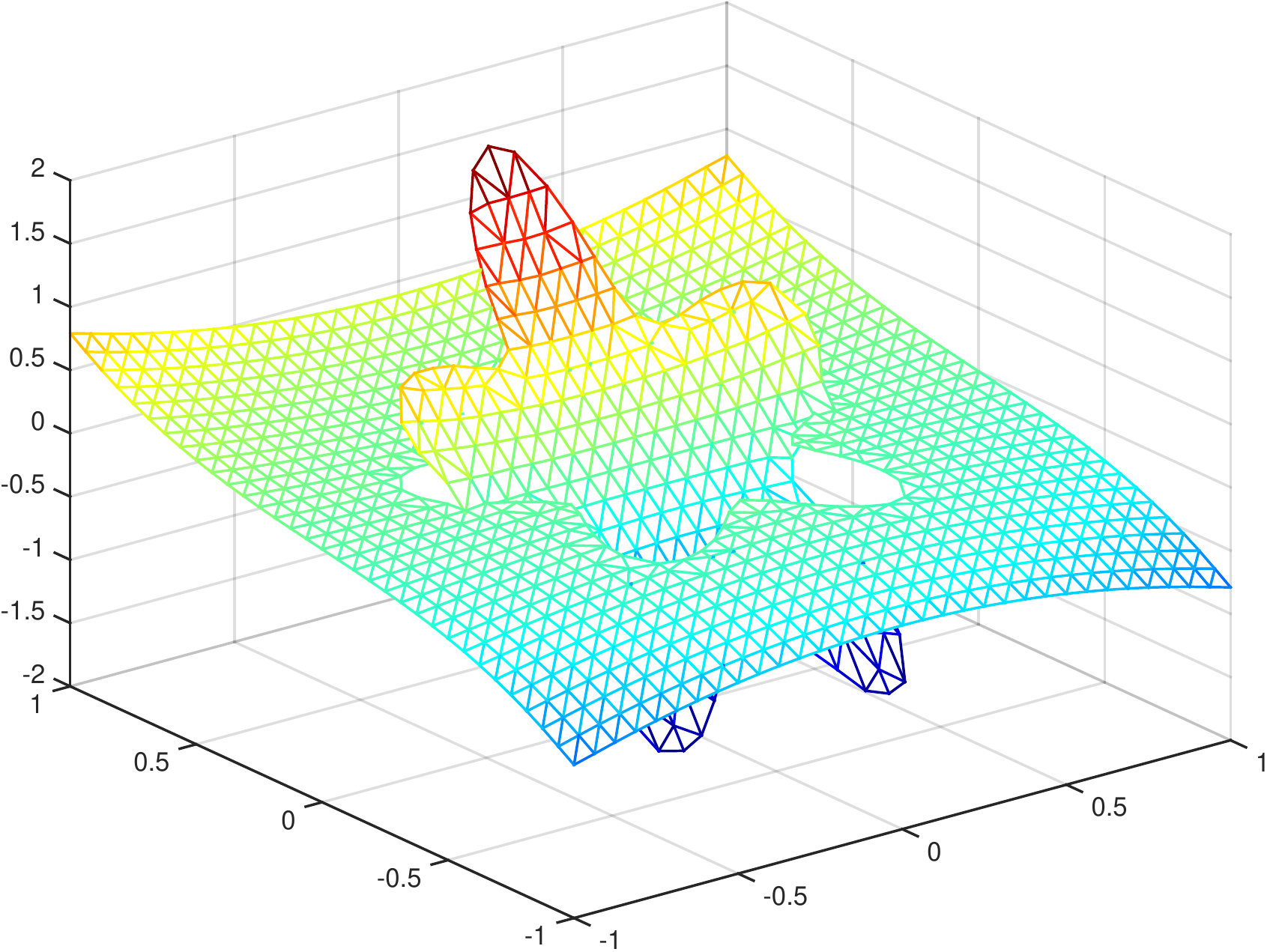}
  \label{fig:flowerry}}
\caption{Plot of  recovered gradient for  Example 5.2: (a) $x$-component; (b) $y$-component.}
\label{fig:flowergradient}
\end{figure}

We use B\"{o}rgers algorithm \cite{Borgers1990} to generate the body-fitted meshes, with the first level of mesh
shown in \cref{fig:flowermesh} and  the finite element solution in \cref{fig:flowersol}.   \cref{fig:flowergradient} gives
the plot of recovered gradient.

In \cref{tab:ex52}, one can see that $De$ decays
at the rate of $\mathcal{O}(h)$, and $D^ie$ superconverges at the rate of $\mathcal{O}(h^{1.5})$ which is consistent with \cref{thm:superclose}.
IPPR ($D^re$) has an $\mathcal{O}(h^{1.5})$ superconvergence that agrees with \cref{thm:superconvergence}.
However, no convergence is observed for standard PPR gradient recovery ($D^pe$) since the exact solution is only smooth on each subdomain.
\begin{table}[htb!]\label{tab:ex52}
\small
\centering
\caption{Convergence rate for Example 5.2. }
\begin{tabular}{|c|c|c|c|c|c|c|c|c|c|}
\hline
 Dof & $De$ & order& $D^{i}e$ & order& $D^{r}e$ & order& $D^{p}e$ &  order\\ \hline\hline
 1089 &7.25e-02&--&1.29e-02&--&1.15e-02&--&6.40e-01&--\\ \hline
 4225 &3.72e-02&0.49&4.67e-03&0.75&3.74e-03&0.83&6.35e-01&0.01\\ \hline
 16641 &1.87e-02&0.50&1.68e-03&0.74&1.19e-03&0.83&6.36e-01&-0.00\\ \hline
 66049 &9.42e-03&0.50&6.06e-04&0.74&3.75e-04&0.84&6.34e-01&0.00\\ \hline
 263169 &4.72e-03&0.50&2.18e-04&0.74&1.27e-04&0.78&6.34e-01&0.00\\ \hline
 1050625 &2.36e-03&0.50&7.69e-05&0.75&4.49e-05&0.75&6.33e-01&0.00\\ \hline
\end{tabular}
\end{table}

{\bf Example 5.3}  This is the same example as used in \cite{Du2016}.  We decompose the computational domain
 $\Omega = (0, 1)\times (0, 1)$ into two parts:
 $\Omega^- = \{ z=(x,y)\in \Omega: x> 0, y>0\}$ and $\Omega^+ = \Omega\setminus \Omega^-$.  The diffusion coefficient $\beta$ in \eqref{equ:model} is chosen as
 \begin{equation*}
\beta(z) =
\left\{
\begin{array}{ccc}
   \beta^- &   \text{if } z\in \Omega^-,\\
       1 &  \text{if } z\in \Omega^+,
\end{array}
\right.
\end{equation*}
with $\beta^-$ as a constant.  When $f=0$ in \eqref{equ:model}, the exact solution $u$ in polar coordinate is given by
 \begin{equation*}
u(r, \theta) =
\left\{
\begin{array}{llc}
    r^{\mu}\cos(\mu(\theta-\pi/4))&  \text{if } 0 \le \theta \le \pi/2, \\
       r^{\mu}\nu\cos(\mu(\theta-5\pi/4))  &   \text{if } \pi/2 \le \theta \le 2\pi,\\
   \end{array}
\right.
\end{equation*}
where
 \begin{equation*}
\mu = \frac{4}{\pi}\left( \sqrt{\frac{3+\beta^-}{1+3\beta^-}} \right) \quad \text{        and        } \quad
\nu = - \beta^-\frac{\sin(\mu\pi/4)}{\sin(3\mu\pi/4)}.
\end{equation*}
Note that $u \in H^{1+s}(\Omega^\pm)$ for any $0 < s < \mu$.

\begin{figure}[ht]
\centering
\subfigure[]{%
     \includegraphics[width=0.4\textwidth]{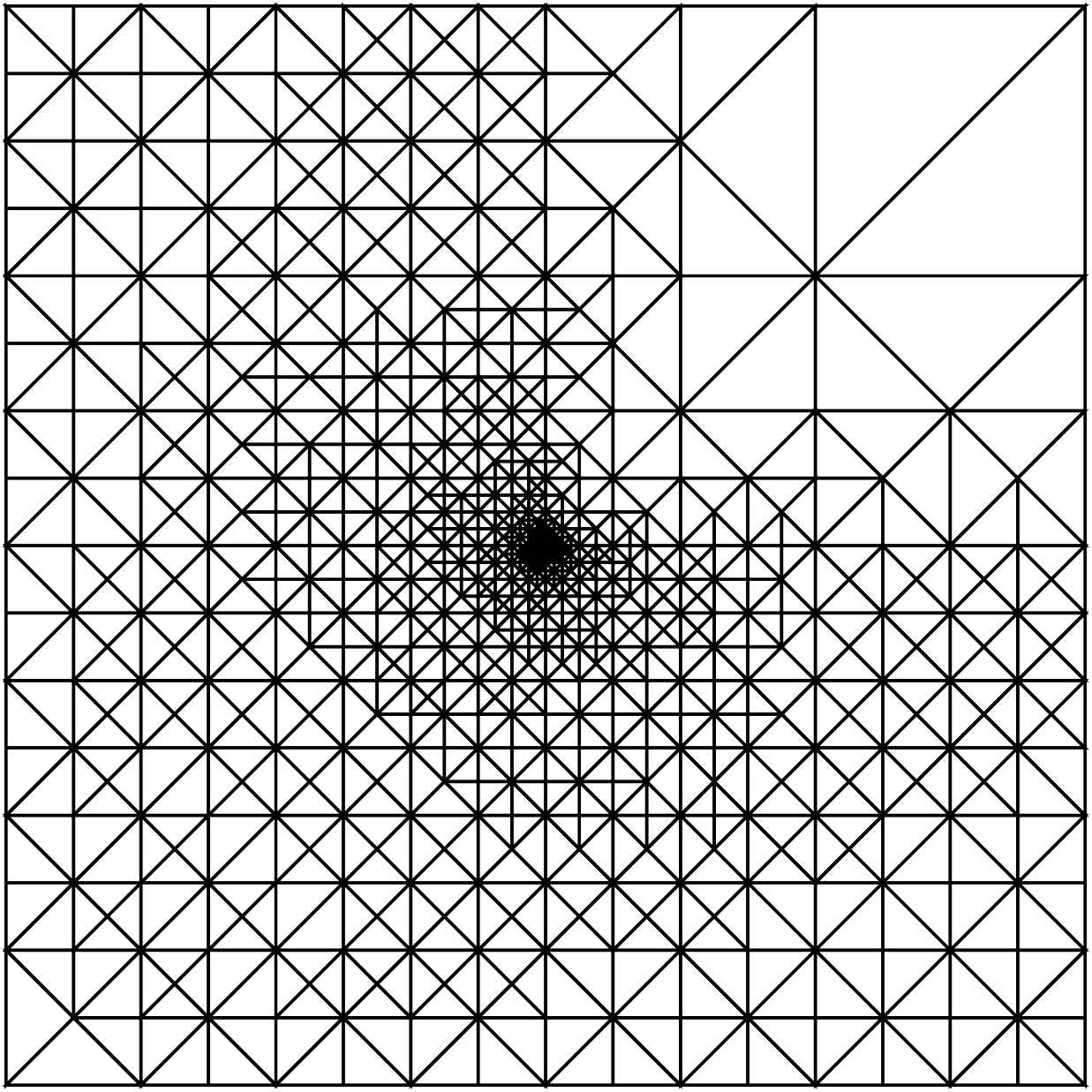}
  \label{fig:Dumesh}}
\quad
\subfigure[]{%
     \includegraphics[width=0.54\textwidth]{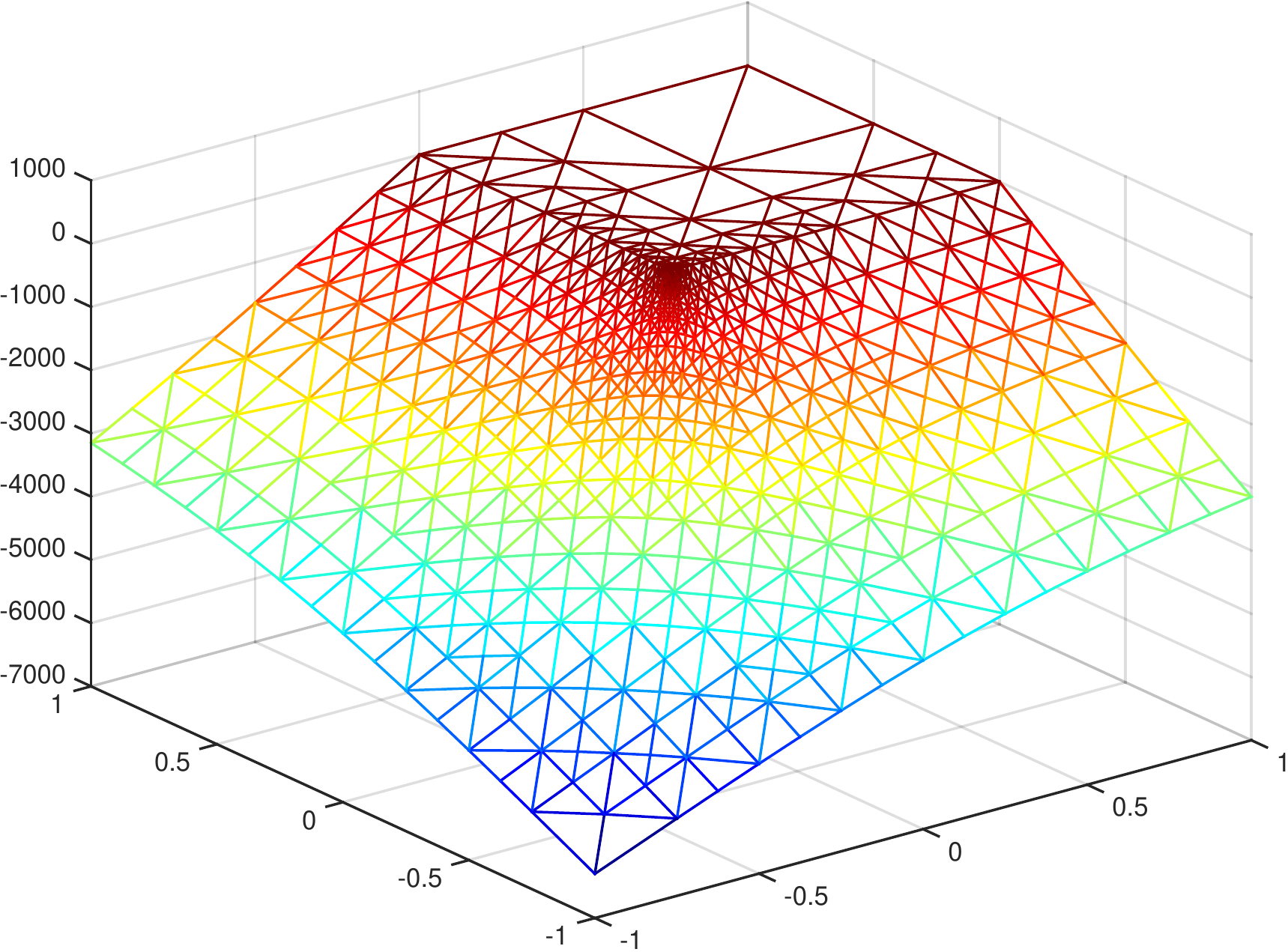}
  \label{fig:Dusol}}
\caption{Adaptive refined  and solution of Example 5.3 with $\beta^-=10000$.  (a)  Adaptive refined mesh. (b).  Finite element solution on adaptive refined mesh h.}
\label{fig:Du}
\end{figure}

\begin{figure}[ht]
\centering
\subfigure[]{%
     \includegraphics[width=0.47\textwidth]{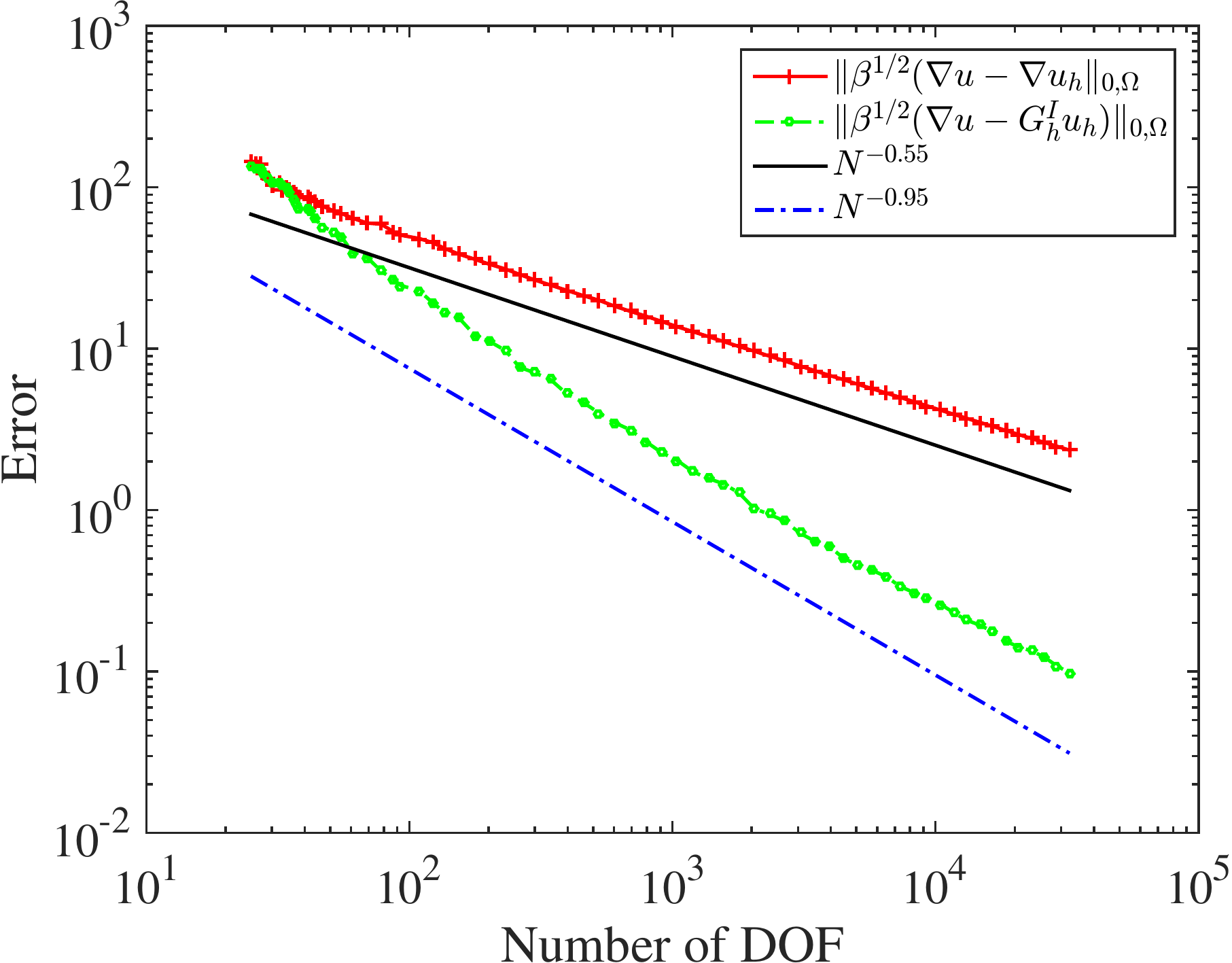}
  \label{fig:Du1000}}
\quad
\subfigure[]{%
     \includegraphics[width=0.47\textwidth]{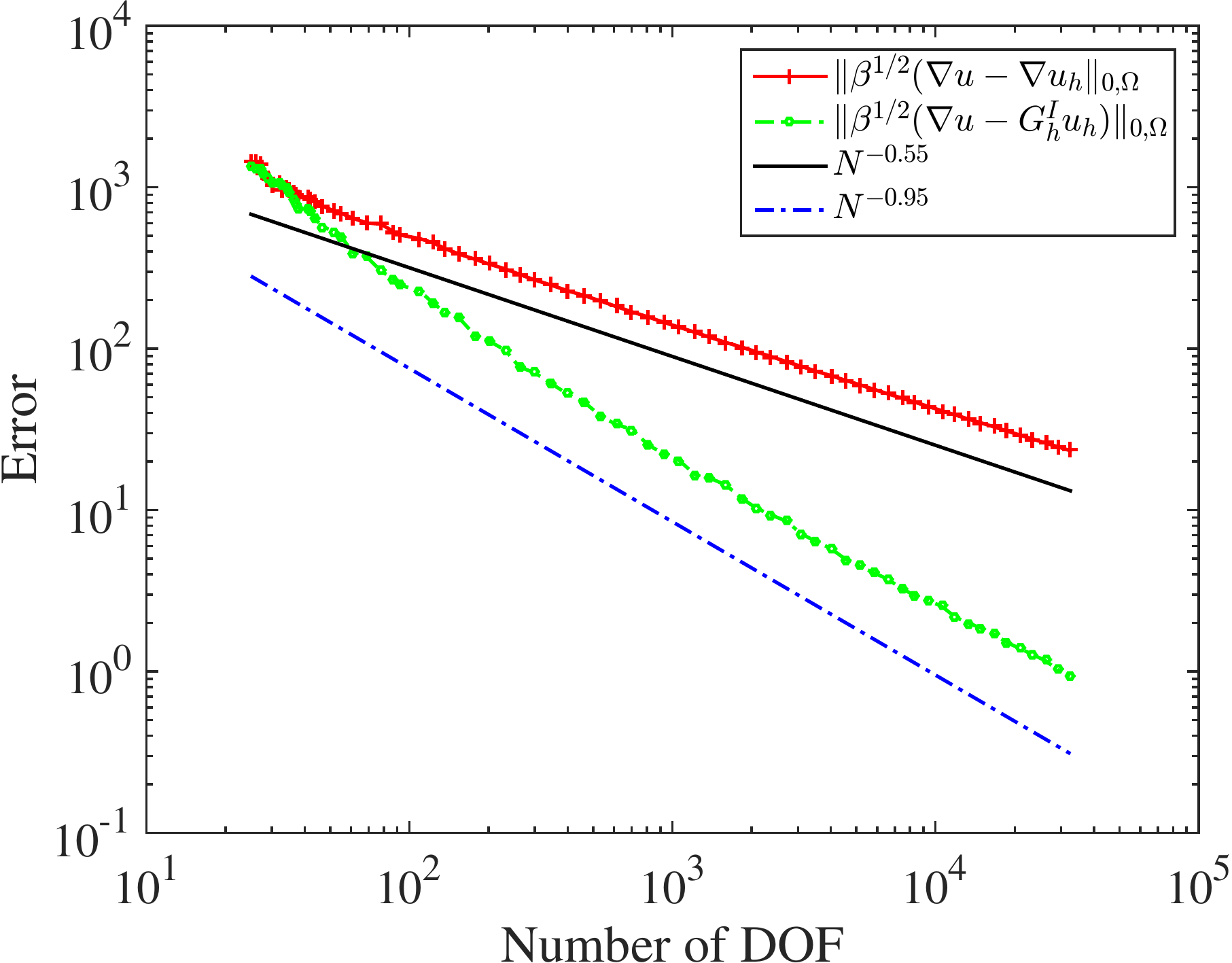}
  \label{fig:Du10000}}
\caption{Convergence rates for Example 5.3:  (a)  $\beta^-=1000$; (b) $\beta^-=10000$.}
\label{fig:Duerr}
\end{figure}

\begin{figure}[!h]
  \centering
  \includegraphics[width=0.7\textwidth]{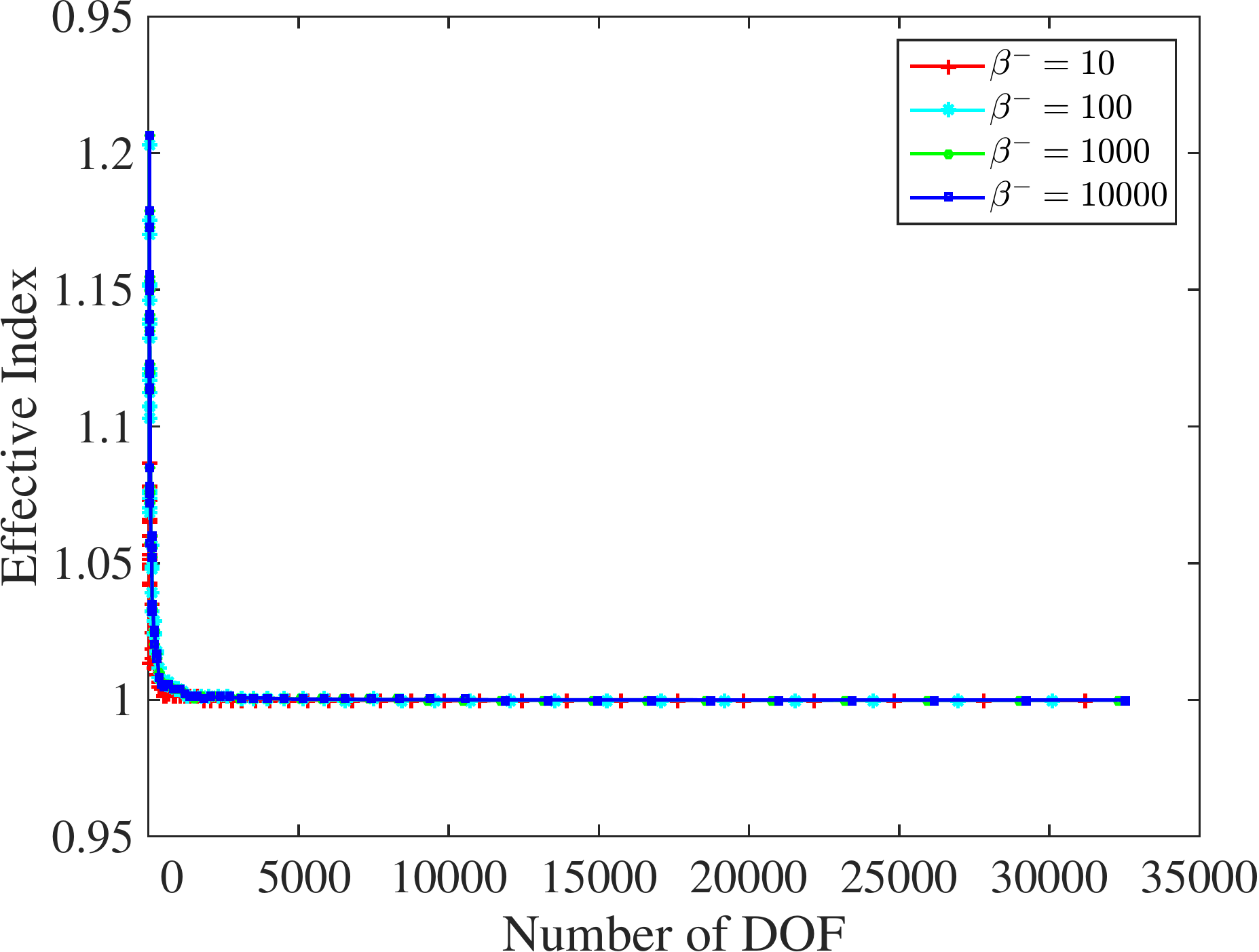}
  \caption{Graph of effective index of Example 5.3}
\label{fig:Du_eff}
\end{figure}

When $\beta^- > 1$, there is a singularity at the origin.  To obtain optimal convergence rate,   we use adaptive finite element method
based on the recovery-type {\it a posteriori} error estimator \cref{equ:localind}.
The bulk marking strategy by \cite{Dorfler1996} with $\theta = 0.2$  is used in numerical computation.
We start with a uniform  initial mesh consisting of $32$ right triangles.
Here, we consider the cases when $\beta^-= 10, 100, 1000$, and $10000$.  \cref{fig:Dumesh} plots  an adaptive refined mesh and
\cref{fig:Dusol} shows the finite element solution when $\beta^-=10000$.   It shows clearly that the refinement is concentrated on  the singularity point.

\cref{fig:Du1000,fig:Du10000} give the numerical convergence rates for $\beta^-=1000$ and $\beta^-=10000$ respectively.
%The numerical results for other value of $\beta^-$ are similar.
In both cases, optimal convergence of $\mathcal{O}(N^{-0.5})$ for energy error  and superconvergence
rate of $\mathcal{O}(N^{0.95})$  can be observed, which is consistent with \cref{thm:superafem}.  We also plot the effective index
\begin{equation*}
\kappa = \frac{\eta_h}{ \|\beta^{1/2}(\nabla u-\nabla u_h)\|_{0, \Omega}}
\end{equation*}
in \cref{fig:Du_eff}.   It shows the error indicator \cref{equ:localind} is an asymptotically exact {\it a posteriori} error estimator for interface problem as in \cref{thm:asy}

\begin{figure}[ht]
\centering
\subfigure[]{%
     \includegraphics[width=0.4\textwidth]{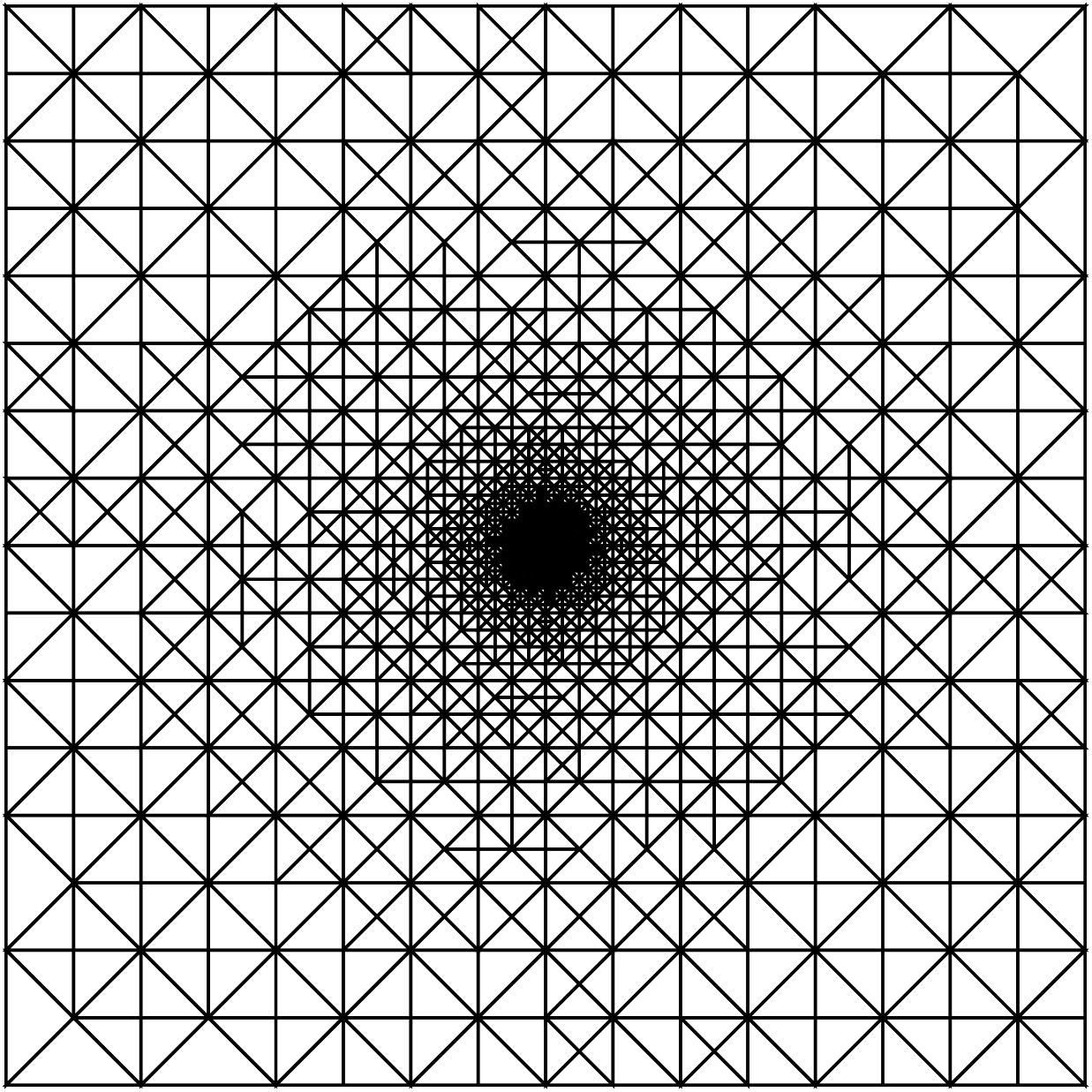}
  \label{fig:kelloggmesh}}
\quad
\subfigure[]{%
     \includegraphics[width=0.54\textwidth]{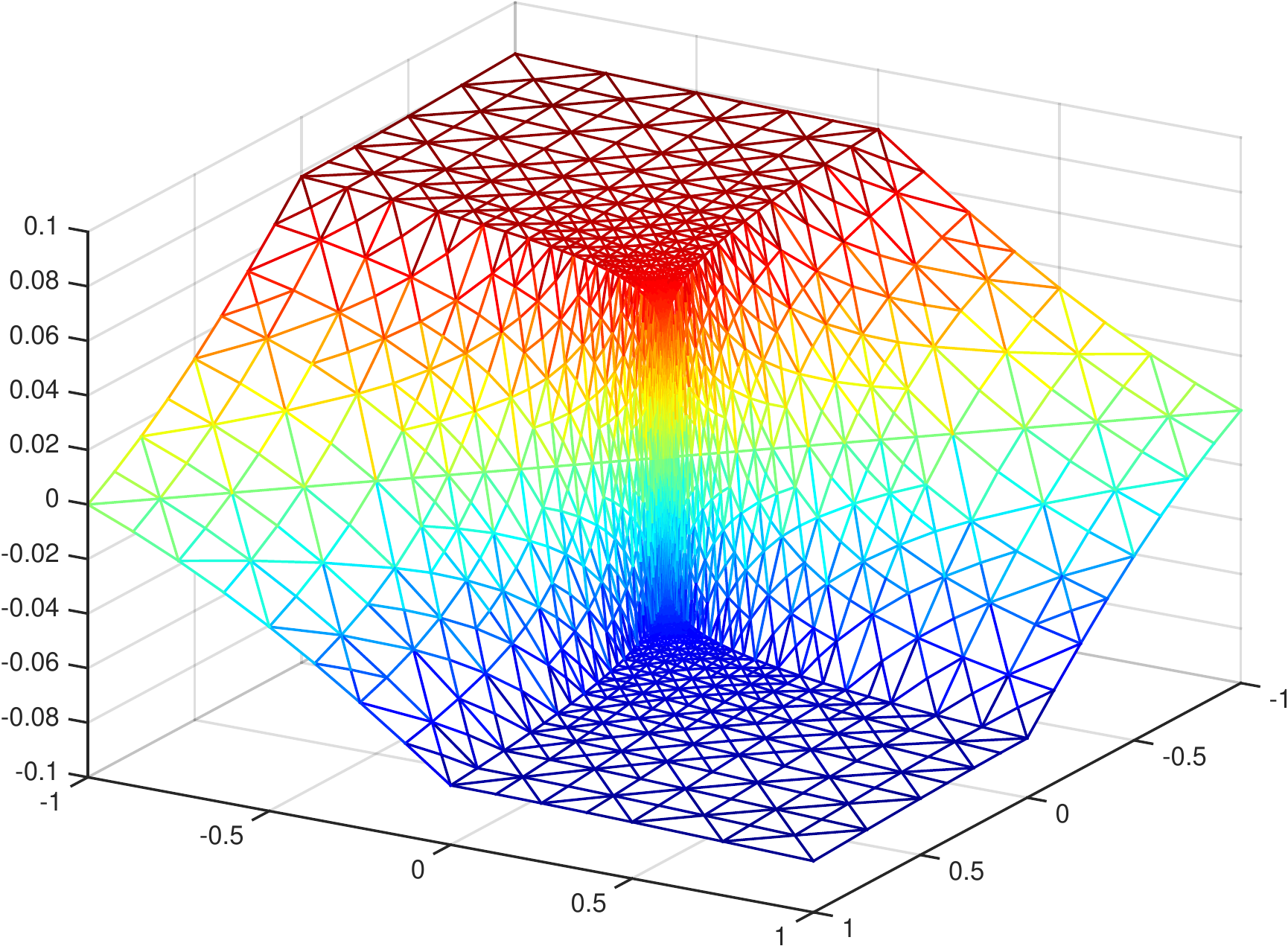}
  \label{fig:kelloggsol}}
\caption{Adaptive refined  and solution of Example 5.3 with $\beta^-=10000$.  (a)  Adaptive refined mesh. (b).  Finite element solution on adaptive refined mesh h.}
\label{fig:kellogg}
\end{figure}

\begin{figure}[ht]
\centering
\subfigure[]{%
     \includegraphics[width=0.47\textwidth]{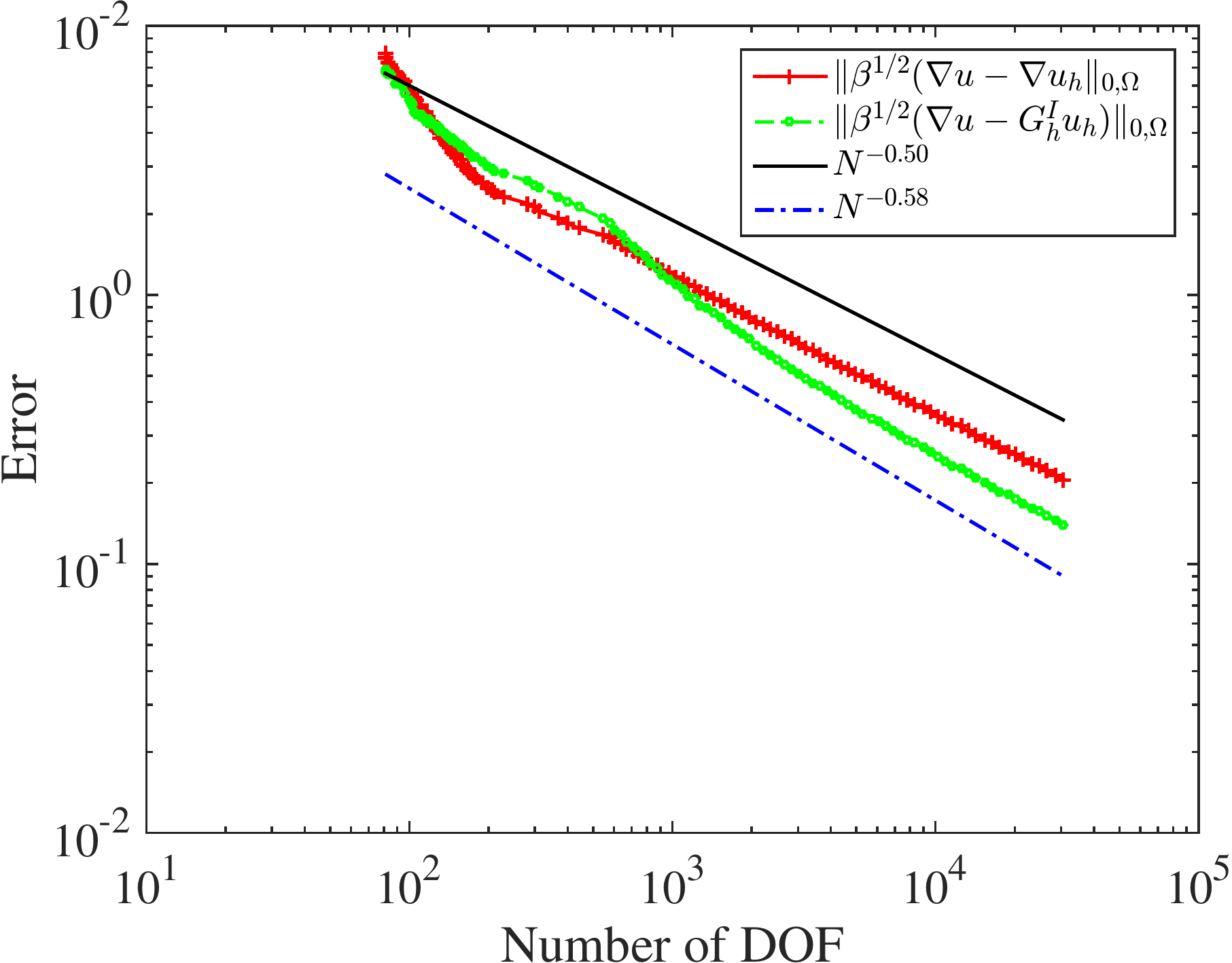}
  \label{fig:kelloggerr}}
\quad
\subfigure[]{%
     \includegraphics[width=0.47\textwidth]{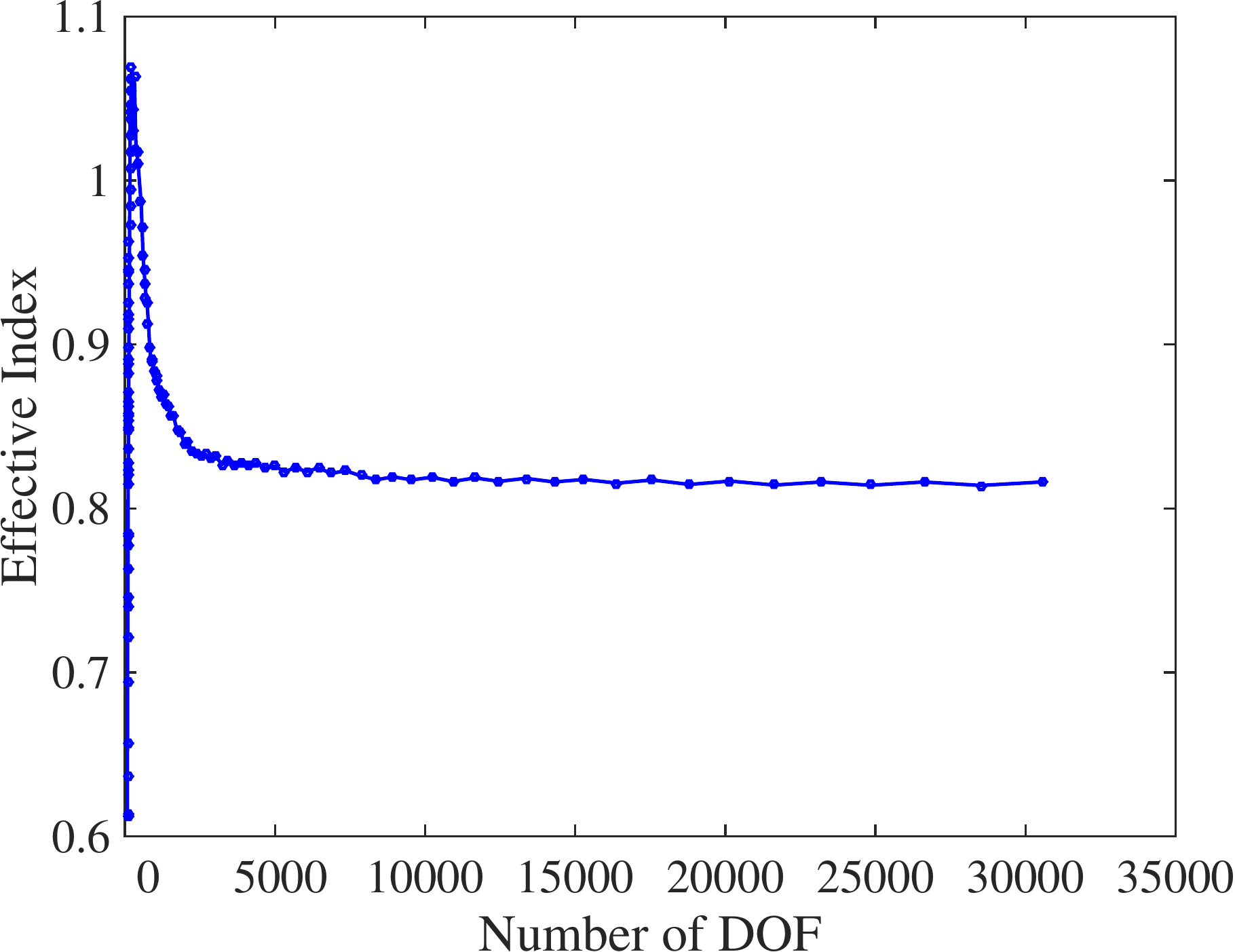}
  \label{fig:kelloggind}}
\caption{Convergence rates for Example 5.3:  (a)  Errors; (b) Effective index.}
\label{fig:kellogghist}
\end{figure}

{\bf Example 5.4.} In the example,   we consider the Kellogg problem which is the benchmark problem of
adaptive finite element method for interface  problem studied by, for example, \cite{CaiZhang2009,ChenDai2002,MorinNochetto2003, Kellogg1974}.
We choose the computational domain $\Omega = (0, 1)\times (0, 1)$, and consider \cref{equ:model}
with $\beta(x) = R$ in the first and third quadrants and $\beta(x) = 1$ in the second and fourth quadrants.
When $f=0$ in \cref{equ:model}, the exact solution $u$ in polar coordinates is given by $u(r, \theta) = r^{\epsilon}\mu(\theta)$  with
\begin{equation*}
\mu(\theta) =
\left\{
\begin{array}{llc}
    \cos((\pi/2-\xi)\epsilon)\cdot \cos((\theta -\pi/2+\nu)\epsilon) &  \text{if } 0 \le \theta \le \pi/2, \\
   \cos(\nu\epsilon)\cdot\cos((\theta-\pi+\xi)\epsilon)  &   \text{if } \pi/2 \le \theta \le \pi,\\
   \cos(\xi\epsilon)\cdot\cos((\theta-\pi-\nu)\epsilon) &  \text{if } \pi \le \theta \le 3\pi/2,\\
   \cos((\pi/2-\nu)\epsilon)\cdot\cos((\theta-3\pi/2-\xi)\epsilon)& \text{if } 3\pi/2 \le \theta \le 2\pi,
\end{array}
\right.
\end{equation*}
with the constants $\epsilon, R, \xi$ and $\nu$ satisfying the nonlinear relations in \cite{ChenDai2002, Kellogg1974}.
Here we choose
\begin{equation*}
\epsilon = 0.1, \quad  \nu = \pi/4, \quad  \xi = -14.9225565104455152, \quad R = 161.4476387975881,
\end{equation*}
and then the exact solution $u \in H^{1+\epsilon}$ with a singularity at the origin.

We start with a uniform  initial mesh consisting of $128$ right triangles and adopt bulk marking strategy by \cite{Dorfler1996} with $\theta = 0.2$.   \cref{fig:kelloggmesh} plots one adaptive refined mesh and \cref{fig:kelloggmesh} plots its corresponding finite element solution.
It clearly indicates that recovery type {\it a posteriori} error estimator  \cref{equ:localind} successfully captures the singularity without  introducing
any overrefinement.   However, the  recovery type {\it a posteriori} error estimator based on  classical gradient recovery operators like SPR or PPR have the problem of overfinement as discussed in \cite{CaiZhang2009}.

\cref{fig:kelloggerr} shows the numerical errors.  One can observe the optimal convergence rate $\mathcal{O}(N^{0.5})$ for energy error and
$\mathcal{O}(N^{0.58})$ superconvergence rate  for recovered energy error.  \cref{fig:kelloggind} gives the history of effective index.
Due to extreme low global regularity of exact solution, the recovery type {\it a posteriori} error estimator  \cref{equ:localind} is not asymptotically exact.
However, it  serves as a robust  {\it a posteriori} error estimator for interface problem as illustrated in \cref{fig:kelloggmesh}.

\section{Conclusion}
In this paper, we develop a novel gradient recovery method for elliptic interface problem based on
body-fitted mesh. Specifically, we define an immersed gradient recovery operator, which overcomes the drawback that stand gradient recovery method fails to produce superconvergence results when solution is lack of regularity at interface. The superconvergence of this method is proved for both mildly unstructured mesh and adaptive mesh.
Several two-dimensional numerical examples are given to confirm our theoretical results, and
verify the robustness of the method served as {\it a posteriori} error estimator. As a continuous study, we plan to develop
gradient recovery methods based on unfitted mesh for elliptic interface problem.

\section*{Acknowledgement}
This work was partially supported by the NSF grant DMS-1418936, KI-Net NSF RNMS
grant 1107291, and Hellman Family Foundation Faculty Fellowship, UC Santa Barbara.

\bibliographystyle{siamplain}
\bibliography{references}
\end{document}